\documentclass[10pt, a4paper]{amsart}

\usepackage{amssymb}
\usepackage{amsmath}
\usepackage{amsthm}
\usepackage{pb-diagram}
\usepackage{enumerate}
\usepackage[latin1]{inputenc}

\theoremstyle{plain} \newtheorem{theorem}{Theorem}[section]
\theoremstyle{plain} \newtheorem{corollary}[theorem]{Corollary}
\theoremstyle{plain} \newtheorem{proposition}[theorem]{Proposition}
\theoremstyle{plain}\newtheorem{lemma}[theorem]{Lemma}
\theoremstyle{definition} \newtheorem{definition}[theorem]{Definition}
\theoremstyle{definition}
\newtheorem{example}[theorem]{Example}
\theoremstyle{remark}
\newtheorem{remark}[theorem]{Remark}
\theoremstyle{remark}
\newtheorem{remarks}[theorem]{Remarks}
\theoremstyle{plain} \newtheorem{conjecture}[theorem]{Conjecture}

\newcommand{\mc}[1]{\mathcal{#1}}

\newcommand{\mbb}[1]{\mathbb{#1}}
\newcommand{\mf}[1]{\mathfrak{#1}}
\newcommand{\C}{\mathcal{C}}
\newcommand{\K}{\mathcal{K}}
\newcommand{\A}[1]{\operatorname{Alb}(#1)}

\begin{document}

\title[Deformations of K\"{a}hler manifolds with vector fields]
{Deformations of K\"{a}hler manifolds with non vanishing holomorphic vector fields}

\author{Jaume Amor\'{o}s, M\`{o}nica Manjar\'{\i}n, Marcel Nicolau}

\address{Jaume Amor\'{o}s. Departament de Matem\`{a}tica Aplicada I,  Universitat Polit\`{e}cnica de Catalunya,
Diagonal 647, E-08028 Barcelona,
Spain} 
\address{M\`{o}nica Manjar\'{\i}n. Institut de Recherche Mathématique de Rennes (IRMAR), 
Campus de Beaulieu Bât.22, 35042 Rennes cedex, France} 
\address{Marcel Nicolau. Departament de
Matem\`{a}tiques, Universitat Aut\`{o}noma de Barcelona, Bellaterra 08193,
Spain}

\medskip

\email{jaume.amoros@upc.edu, monica.manjarin-arcas@univ-rennes1.fr,\newline
nicolau@mat.uab.cat}

\thanks{This work was partially supported by the Ministerio de Educaci\' on y Ciencia of Spain,
grants MTM2008-02294 and MTM2009-14163-C02-02, and by the Generalitat de Catalunya grants
2009SGR 1207 and 2009SGR 1284. The second author was supported by the Programa de Movilidad de 
Recursos Humanos del Plan Nacional de I-D+I 2008/2011.}



%%% ----------------------------------------------------------------------

\begin{abstract}

We study compact K\"{a}hler manifolds $X$ admitting nonvanishing
holomorphic vector fields, extending the classical birational classification of 
projective varieties with tangent vector fields to a classification modulo deformation
in the K\"ahler case, and biholomorphic in the projective case. 
We introduce and analyze a new class of {\em tangential deformations}, and show that 
they form a smooth 
subspace in the Kuranishi space of deformations of the complex structure of $X$. 
We extend Calabi's theorem on 
the structure of compact K\"{a}hler manifolds $X$ with $c_1(X) =0$
to compact K\"ahler manifolds with nonvanishing tangent fields,
proving that any such manifold $X$
admits an arbitrarily small tangential deformation which is a 
suspension over a torus; that is,
a quotient of $F\times \mbb C^s$ fibering over a torus $T=\mbb C^s/\Lambda$. 
We further show that either $X$ is uniruled or, up to a finite Abelian covering, 
it is a small deformation of a product $F\times T$
where $F$ is a K\"{a}hler manifold without tangent vector fields and $T$ is a torus. 
A complete classification when $X$ is a projective manifold, in which case the deformations
may be omitted, or when $\dim X\leq s+2$ is also given.
As an application, it is shown that the study of the dynamics of holomorphic
tangent fields on compact K\"ahler manifolds reduces to the case of rational
varieties. 
\end{abstract}

%%% ----------------------------------------------------------------------
\maketitle
%%% ----------------------------------------------------------------------

\tableofcontents

\section*{Introduction}

This article is devoted to the study of compact K\"{a}hler manifolds admitting nonvanishing
holomorphic vector fields or, what is equivalent, endowed with a locally free holomorphic $\mbb C^s$-action
for a certain $s>0$. 
We introduce a particular type of deformations of such manifolds and, 
as an application, we describe quite precisely their structure. 
Our study relies on the work carried out simultaneously and independently by A.~Fujiki \cite{Fuji2} and 
D.~Lieberman \cite{Lieb}. We use specially their deep result on the 
structure of the group of holomorphic automorphisms of a compact K\"{a}hler manifold 
(cf. Theorem \ref{thm:fuji-lieb}).
\smallskip

Our work is a continuation of the classical study of complex projective 
manifolds with tangent vector fields, which led to their birational
classification up to finite covering, successively developed by F. Severi, R. Hall 
and D. Lieberman (\cite{Severi}, \cite{Hall}, \cite{Lieberman}).
The conclusions were summarized in Theorem~1 of \cite{Lieberman} as

\begin{theorem} \label{thm:classic}
Let $X$ be a complex projective manifold with a holomorphic tangent vector field $v$.
Then:
\begin{enumerate}[\rm a)]
\item $X$ has a finite \'etale Abelian covering $X'$ which is birational to a product 
$M \times T \times \mbb C P^r$, where $M$ is another projective manifold, $T$ an Abelian
variety, and $\mbb C P^r$ the projective space of dimension $r \ge 0$, such that the
lift of the tangent vector field $v$ to $X'$ has trivial component in $TM$.
\item If the tangent vector field $v$ on $X$ vanishes at some point then $X$ is ruled.
\end{enumerate}
\end{theorem}

In Theorem \ref{thm:susp_proj} below we continue this classification in the case of 
nonvanishing vector 
fields, obtaining a classification up to biholomorphism by using suspensions
over Abelian varieties (see Example \ref{example_suspension}) and repeated
\'etale Abelian coverings in order to exhaust all tangent vector fields.

Tangent vector fields play a role in the classification of algebraic varieties, 
specially in the case of Kodaira dimension 0. The earliest indication of this fact
known to the authors is the result stated by Eugenio Calabi in~\cite{Cal}.
In that article, E. Calabi formulated his celebrated conjecture about
K\"ahler manifolds admitting Ricci--flat metrics and, under the assumption that 
the conjecture was true, he proved the result stated as Theorem~\ref{thm:calabi2} below.
This theorem describes the structure of compact K\"{a}hler manifolds with trivial 
first Chern class, i.e.  with $c_1(X)=0$. It  was proved for projective manifolds 
by Y.~Matsushima in \cite{Matsu} using previous results by A.~Lichnerowicz, 
but without the assumption of Calabi's conjecture, and extended to the general K\"{a}hler 
case by F.A.~Bogomolov~\cite{Bogo} and independently by 
Lieberman \cite{Lieb}.  

\begin{theorem} [Calabi, \cite{Cal}, \cite{Matsu}, \cite{Bogo}, \cite{Lieb}] \label{thm:calabi2}
Let $X$ be a compact K\"ahler manifold with $c_1(X)=0$ and $b_1(X) >0$. 
Then $X$ admits, as a finite \'etale covering, a product $X'= F \times T$, 
with $T$ a complex torus of real dimension $b_1(X')$ and $F$ a K\"ahler manifold 
with $c_1(F)=0$ and $b_1(F)=0$. Moreover, $X'$ is a regular covering space of 
$X$ and the group of covering transformations is solvable. 
\end{theorem}

This Theorem was improved by Beauville in~\cite{Bea2}, where he completely describes the topology
and the geometry of the fiber $F$. 

The proof of Calabi's theorem may be succintly described as follows: Use the 
triviality of the canonical bundle to show that $X$ (or a finite covering of it) has 
$q=\frac{1}{2} b_1(X)$ linearly
independent nonvanishing tangent vector fields; check that the Albanese morphism is 
onto, the tangent vector fields are lifts from the tangent bundle on the Albanese torus,
and they make $X$ (or a finite covering of it) a suspension with trivial monodromy over 
the Albanese torus. 

In this work we develop a refinement of the above argument to prove an extension to 
general K\"ahler manifolds of Calabi's theorem, namely our Theorem~\ref{thm:susp_proj} and
Proposition~\ref{prop:dicotomia} may be condensed as:

\begin{theorem} \label{teo:mainmain}
Let $X$ be a compact K\"ahler manifold with $s$ non vanishing
tangent vector fields $v_1, \dots, v_s$, defining a locally free $\mbb C^s$-action. 
Then $X$ admits a small deformation 
$X_{\epsilon}$, which is a suspension
over an $s$--dimensional torus $T$.

Furthermore, there is a finite Abelian covering $X'_\epsilon$ of a small deformation 
$X_\epsilon$ of $X$, 
a torus $T'$ of dimension $s'\geq s$ and a compact K\"{a}hler manifold $F$ without 
non vanishing vector fields in such a way that, \begin{enumerate}[\rm (i)]
\item if $\operatorname{kod}(X)\geq 0$ then $X'_\epsilon = F\times T'$ and $F$ has no vector fields,
\item if $\operatorname{kod}(X)=-\infty$ then $X'_\epsilon$ is a topologically trivial suspension over $T'$ and fiber $F$.
\end{enumerate}
In both cases $\operatorname{kod}(F)=\operatorname{kod}(X)$. If $X$ is a complex projective manifold one may omit the small deformation of the complex structure, and the above
decomposition holds for $X_\epsilon=X$. 
\end{theorem}

Here, $\operatorname{kod}(X)$ denotes the Kodaira dimension of $X$. The definition of a suspension 
manifold is recalled in Example~\ref{example_suspension}.
\smallskip

We notice that a locally free holomorphic $\mbb C^s$-action, i.e. a holomorphic $\mbb C^s$-action 
for which every isotropy group is discrete, is determined by $s$ holomorphic vector fields 
$v_1, \dots, v_s$ generating 
an $s$-dimensional Abelian Lie algebra $\mf a$ such that each vector field $v$ in $\mf a$
different from zero is a non vanishing vector field. On a compact Kahler manifold $X$, 
the set $\mf h^1$ of vector fields having zeros is an ideal of the Lie algebra $\mf h$ 
of holomorphic vector fileds on $X$, there is a direct sum decomposition  $\mf h = \mf h^1 \oplus \mf a$,
where $\mf a$ is an Abelian subalgebra of $\mf h$, and $X$ admits a locally free 
holomorphic $\mbb C^s$-action if and only if $\dim \mf a >0$ 
(cf. Propositions~\ref{prop:LieSubA} and \ref{prop:touzet}). 
Hence the hypothesis in the above Theorem is fulfilled 
whenever $s = \dim (\mf h/\mf h^1) > 0$. This is the case, in particular, if $X$ admits a 
non vanishing vector field.

\smallskip

The small deformation of the complex structure of $X$ cannot be avoided in the general K\"ahler case 
if $c_1(X) \neq 0$, as it is shown by the examples discussed in Section~\ref{sec:structure}.

On the other hand, for $X$ compact K\"ahler with $c_1(X)=0$, the dimension of the space of holomorphic
vector fields on $X$ is $b_1(X)/2$ and this fact also renders innecessary the small deformation 
of complex structure, thus Corollary~\ref{cor:c1=0} to our Theorem~\ref{teo:mainmain} 
refines Calabi's theorem, as we show that the group of covering transformations is in fact Abelian. 
This refinement had been previously established in \cite{Bogo} in the case when $X$ is a
complex projective manifold.

\smallskip

Besides this extension of Calabi's theorem from Calabi--Yau manifolds to 
any K\"ahler manifolds with nonvanishing tangent vector fields, the authors 
hope that the results of this
work may have further applications in the case of manifolds of Kodaira dimension zero.
K. Ueno (cf.~\cite[\S11]{Ue}) made a conjecture, which is the extension of Calabi's one to this setting, 
and was refined 
by J\'anos Koll\'ar (cf. {\cite{Kol} 4.16}) as:

\begin{conjecture} [Ueno, Koll\'ar] \label{conj:kollar}
Let $X$ be a smooth, projective variety with $\text{kod}(X)=0$.
Then $X$ has a finite, \'etale covering $\tilde X$ such that $\tilde X$
is birational to the product of an Abelian variety and a simply
connected variety $F$ with $\text{kod}(F)=0$.
\end{conjecture}

Corollary \ref{c:k=0} is a modest contribution to this problem, establishing the conjecture
if one has $\dim X -2$ linearly independent nonvanishing tangent vector fields, a particular 
case being threefolds with non vanishing vector fields. 

This Corollary is a particular case of the classification of compact K\"{a}hler and projective manifolds 
with $\dim X -2$ linearly independent nonvanishing tangent vector fields carried out in Section~\ref{s:small_codimension} based on Theorem~ \ref{teo:mainmain} and the classification 
of curves and surfaces.
\smallskip

In Section \ref{s:dinamica} we extend the results in \cite{Lieb} concerning 
the dynamical properties of a holomorphic tangent vector field $v$ on a 
compact K\"ahler manifold $X$. Applying Theorem~\ref{teo:mainmain},
we show that the dynamical system $(X,v)$ becomes trivial after taking a finite unramified covering
space of $X$ if $\operatorname{kod}(X) \ge 0$ and that, in all cases, the dynamics of the system
reduce to the dynamics of an Abelian group $\mbb C^p\times(\mbb C^*)^q$ acting on a 
rational variety (Theorem~\ref{teo:dinamica_racional}).

\bigskip

In order to achieve all the results described above, we
introduce the notion of tangential deformation of a holomorphic group action.
A locally free holomorphic action of a (connected) complex Lie group $G$ on a complex manifold 
$X$ defines a  holomorphic foliation $\mc F$ on $X$ whose leaves are the orbits of the action. 
$G$-equivariant deformations of $X$ change the complex structure of $X$ as well as the holomorphic 
foliation $\mc F$.
The existence of a versal space of deformations for equivariant deformations of compact complex manifolds 
endowed with a holomorphic action was proved by Cathelineau in~\cite{Cat}.
In this article we consider equivariant deformations of the locally free $G$-action for which the 
foliation $\mc F$ keeps fixed its holomorphic transverse structure. This type of deformations, called
{\sl tangential} deformations of the action, also admits a versal space of deformations whose tangent 
space at the origin is naturally identified to the space of infinitesimal tangential deformations
(Theorem~\ref{teo:kuranishi_equivariant}).

We focus on the case in which $X$ is a K\"{a}hler manifold. Under that hypothesis and as a consequence of 
the mentioned theorem by Fujiki and Lieberman, the group $G$ is necessarily Abelian. 
For a given locally free holomorphic $\mbb C^s$-action on a compact K\"{a}hler manifold $X$ we construct 
a family $X_R$ of tangential deformations of the action parametrized by smooth space 
$R$ with the following properties:
\begin{enumerate}[(1)]
\item $R$ is an open subset of an Euclidean space $\mbb C^N$, complementary to an affine real
algebraic variety.
\item The family $X_R$ is versal at each point $r\in R$ (Theorem~\ref{teo:versal_space} and
Proposition~\ref{prop:tot_versal}).
\item Each element $X_r$, with $r\in R$, is a K\"{a}hler manifold (Theorem~\ref{teo:kahler_stability}).
\end{enumerate}
We emphasize that $R$ is a (real) Zariski open set and not merely an small open ball 
of $\mbb C^N$.

Looking at general properties of tangential deformations of the action, we prove that the 
Kodaira dimension $\operatorname{kod}(X)$ of $X$ is constant under tangential deformations.
We also relate $R$ to the Kuranishi space $\K_X$ of $X$ showing that the forgetful map
$$
R\longrightarrow \K_X
$$
has a smooth image of dimension $s\cdot b_1(X)/2$. In particular this shows the following 
(cf. Theorem~\ref{teo:kuranishi})

\begin{theorem}
The space $H^1(X, \Theta_X)$ of infinitesimal deformations of a compact K\"{a}hler manifold $X$
endowed with a locally free  $\mbb C^s$-action contains a subspace of dimension 
$\ell= s\cdot b_1(X) /2$ of unobstructed infinitesimal deformations. Consequently the Kuranishi space
$\K_X$ contains a smooth subspace of dimension $\ell$.
\end{theorem}

Many of the results of this article are also valid for manifolds $X$ belonging to the 
class $\mc C$, the class introduced by A.~Fujiki and whose elements are those manifolds that 
are bimeromorphic to K\"{a}hler manifolds. Hence some of the results are stated in that context. 
However, we do not dispose of a characterization of those suspension manifolds 
that are in the class $\mc C$ (cf. Remark~\ref{rem:prob_fujiki}). 
Moreover, the class $\mc C$ is not stable under small deformations 
and by these reasons we do not know if the structure theorems of Section~\ref{sec:structure} 
are valid for Fujiki manifolds and we are able to prove them only in the case of K\"{a}hler manifolds.

Even if it is not said explicitly, all the group actions considered along this article are 
holomorphic and the covering spaces are always unramified (\'{e}tale).

We are very indebted to M. Brunella for inspiring discussions and to F. Touzet who
communicated the proof of Proposition~\ref{prop:touzet} to us.

%%%%%%% SECTION 1

\section{Groups of automorphisms of K\"{a}hler manifolds} \label{Sec:Prel}

Let $X$ be a compact complex manifold of dimension $n$. All along this article, 
with the only exception of Section~\ref{Sec:TanDef},
the manifold $X$ will be assumed to be K\"{a}hler or, more generally, to belong 
to the class $\C$ introduced by A. Fujiki in \cite{Fuji1}. We recall that  a manifold 
belongs to the class $\C$ (or Fujiki class) if it is bimeromorphic to a K\"{a}hler 
manifold. Notice that 
Moishezon manifolds are in $\C$ and that a submanifold of a Fujiki manifold is 
again in the class $\C$.
Recently, Demailly and Paun have characterized compact manifolds in 
$\C$ as those manifolds carrying a K\"{a}hler current (cf. \cite{DePa}). 

Compact manifolds $X$ in the class $\C$ fulfill many of the properties of compact K\"{a}hler 
manifolds. In this article we will make use of the following ones: (i) holomorphic forms on $X$
are closed, (ii)  there are $\mbb C$-antilinear isomorphisms $H^q(X,\Omega_X^p) 
\equiv H^p(X,\Omega_X^q)$,  
(iii) the Hodge theorem holds, i.e. $H^{m}(X,\mbb C)\equiv\oplus_{p+q=m}H^q(X,\Omega_X^p)$ 
and (iv)  $H^1(X, \mbb C)$ is isomorphic to  $H^1(\A{X}, \mbb C)$, where $\A{X}$ is 
the Albanese torus of $X$. For the proof of the above statements and a general account on 
the properties of Fujiki manifolds we refer to \cite{Fuji2} and \cite{Fuji3}.
In a sharp contrast, the class $\C$ is not stable under small deformations (cf. \cite{Cam} and \cite{LePo}).

Let us denote by
$\mf h_X$ the Lie algebra of holomorphic vector fields on $X$ and let $\mf h_X^{1}$
be the subspace of $\mf h_X$ whose elements are the vector fields on $X$
annihilated by all the holomorphic 1-forms on $X$.
Since global holomorphic 1-forms are closed, $\mf h_X^{1}$ is an ideal of $\mf h_X$ containing 
the derived subalgebra $[\mf h_X,\mf h_X]$ of $\mf h_X$.

The following characterization of $\mf h_X^{1}$ obtained by Carrell and Lieberman 
in \cite{CaLie} was extended by Fujiki to the class $\mc C$ in \cite{Fuji1} .

\begin{theorem}[Carrell-Lieberman \cite{CaLie}]\label{alpha}
Let $X$ be a compact K\"{a}hler manifold. Then $\mf h_X^{1}$ is the vector space of 
holomorhic vector fields with zeros. In particular, if $X$ admits a non vanishing holomorphic vector 
field $v$, then there is a holomorphic 1-form $\alpha$ on $X$ such that $\alpha(v)=1$.
\end{theorem}

Let $\phi: X\rightarrow \A{X}$ be the natural map from $X$ to its
Albanese torus $\A{X}$. Each element of the group $\operatorname{Aut}_{\mbb{C}}(X)$ of
holomorphic automorphisms of $X$ induces an automorphism of $\A{X}$ and the correspondence
$$
\Phi: \operatorname{Aut}_{\mbb{C}}(X) \longrightarrow \operatorname{Aut}_{\mbb{C}}(\A{X})
$$
is a group morphism. The map $\phi$ is equivariant with respect to the action of
$\operatorname{Aut}_{\mbb{C}}(X)$ on $\A{X}$ induced by $\Phi$. The connected component
of the identity $\operatorname{Aut}^0_{\mbb{C}}(X)$
of $\operatorname{Aut}_{\mbb{C}}(X)$ is mapped by $\Phi$ into
$\operatorname{Aut}^0_{\mbb{C}}(\A{X})\cong \A{X}$ and its image is a
sub-torus $T_X$ of $\A{X}$. Notice that $T_X$ is contained in the image $\phi(X)$ of $X$ 
into $\A{X}$.

The following result was obtained independently and at the same time 
by A.~Fujiki and by D.~Lie\-ber\-man. The latter proved it for K\"ahler manifolds and the former 
for manifolds in the class $\mc C$. 

\begin{theorem} [Fujiki \cite{Fuji2}, Lieberman \cite{Lieb}] \label{thm:fuji-lieb}
Let $X$ be a compact complex manifold that is K\"{a}hler or, more generally, that belongs to the 
Fujiki class. Then there is an exact sequence of groups
\begin{equation}\label{successio-grups}
1 \rightarrow L \rightarrow \operatorname{Aut}^0_{\mbb{C}}(X) \overset{\Phi}\rightarrow T_X
\rightarrow 1
\end{equation}
where $L$ is a linear algebraic group (so with finitely many connected components)
with Lie algebra $\mf h_X^{1}$.
Moreover, there is a group morphism
$\iota :\mbb C^r\rightarrow \operatorname{Aut}^0_{\mbb{C}}(X)$ such that the composition
$\Phi\circ\iota$
is the universal covering of $T_X$.
\end{theorem}

\begin{corollary}
If $\mf h_X^1 = 0$ then $\operatorname{Aut}^0_{\mbb{C}}(X)$ is a complex torus.
\end{corollary}

Notice that Theorem \ref{alpha} is also a corollary of Theorem \ref{thm:fuji-lieb}.
Another consequence of the theorem is the following:

\begin{proposition}\label{prop:LieSubA}
There is an Abelian Lie subalgebra $\mf a$ of $\mf h_X$ such that
\begin{equation}\label{descomposicio}
\mf h_X = \mf h_X^{1} \oplus \mf a,
\end{equation}
where the direct sum is in the sense of vector spaces.
\end{proposition}

\begin{remarks} \label{rem:camps} (a) We denote $s= \dim_{\mbb C}T_X 
= \dim_{\mbb C} \mf a = \dim_{\mbb C} \mf h_X/\mf h_X^{1}$. 
Then we have
$$
s\leq \dim_{\mbb C}\A{X} =b_1(X)/2 \quad \mbox{and} \quad s\leq \dim_{\mbb C}(X).
$$

(b) The subalgebra $\mf a$ is not unique. Such a Lie subalgebra
can be characterized as a vector subspace of  $\mf h_X$ of maximal dimension
generated by non vanishing commuting vector fields.

(c) Let $v_1, \dots, v_s$ be a basis of $\mf a$. Then the vector fields $v_1, \dots, v_s$
are linearly independent at each point of $X$.

(d) The choice  of  a basis $v_1, \dots, v_s$ of $\mf a$ determines an
action of $\mbb C^s$ on $X$ which is locally free, i.e. an action which is injective at 
the Lie algebra level or, what is equivalent, which has discrete  isotropy groups. 
Conversely, given such a locally free action 
$\varpi: \mathbb C^m\times X\rightarrow X$,
of maximal rank (i.e. $m$ is maximal), one has $m=s$ and the fundamental
vector fields of the action generate a subalgebra of
$\mf h_X$ complementary to $\mf h^{1}_X$.
If the orbits are given by a torus action, i.e. if the action
factorizes through a torus $T=\mbb C^s/\Lambda$, then $\mf a$ is contained
in the center of $\mf h_X$.
\end{remarks}

The following result is stated without proof in \cite{Lieb}. The proof we give here
is due to F. Touzet.

\begin{proposition}\label{prop:touzet}
The Abelian subalgebra $\mf a$ of $\mf h_X$ fulfilling~$(\ref{descomposicio})$ 
can be chosen in the center of 
$\mf h_X$.
\end{proposition}

\begin{proof}
First we notice that $L$ can be assumed to be connected just by replacing $T_X$ by a 
finite covering, which is again a torus. Let  $\operatorname{Aut}^0(\mf h^1_X)$
be the identity component of the automorphism group of the complex Lie algebra 
$\mf h^1_X = \operatorname{Lie}(L)$. Since the group $L$ is normal in 
$G = \operatorname{Aut}^0_{\mbb{C}}(X)$, the adjoint representation
$$\begin{array}{rcc}
\psi: G &\rightarrow& \operatorname{Aut}^0(\mf h^1_X) \\
g & \mapsto & \operatorname{Ad}_g
\end{array}$$
is well defined and induces a group morphism
$$
\tilde{\psi}:T_X \longrightarrow H = \frac{\operatorname{Aut}^0(\mf h^1_X)}{\operatorname{Ad}(L)} .
$$
Notice that $\operatorname{Ad}(L)$ is closed in $\operatorname{Aut}^0(\mf h^1_X)$,  
as $L$ is an algebraic group, and that $\operatorname{Ad}(L)$ is a normal subgroup of 
$\operatorname{Aut}^0(\mf h^1_X)$,  as $\operatorname{Lie}(\operatorname{Ad}(L))$ 
is an ideal of 
$\operatorname{Der}^0(\mf h^1_X)$.
Therefore the connected group $\operatorname{Ad}(L)$ is a normal and closed
subgroup of the algebraic group $\operatorname{Aut}^0(\mf h^1_X)$.
Hence the quotient $H$ is also a linear algebraic group (cf. \cite{Hump}) and
the image by $\psi$ of the compact group  $T_X$ is necessarily constant, thus
equal to the identity.
As a consequence, $\psi$ has image contained in 
$\operatorname{Ad} (L)$. We deduce that, for each $v\in \mf h_X$, there is $v_1=v_1(v)\in \mf h^1_X$ such that
$$[v,w]=[v_1,w] \quad \forall w\in \mf h^1_X,$$
that is, $v-v_1$ commutes with all the elements in $\mf h^1_X$.

We consider now $v\in \mf h_X$ with the property that $\Phi(\{\operatorname{exp} tv\})$ is dense in
the torus $T_X$. Then $v-v_1$ has the same property.  Hence the
closure $A_{v-v_1}$ of the group $\{\operatorname{exp} t(v-v_1) \}$ is an Abelian subgroup
of $G$ that is mapped onto $T_X$ by $\Phi$. This implies that $v-v_1$
is in the center of $\mf h_X$ and therefore that $A_{v-v_1}$ is in the center of $G$.
Now we take a group
morphism $\iota :\mbb C^r\rightarrow A_{v-v_1}$ with the property that the composition $\Phi \circ \iota$
is the universal covering of $T_X$. Then the morphism $\iota$ determines an Abelian subalgebra 
$\mf a$ of $\mf h_X$ with the requiered properties.
\end{proof}

\begin{remarks}
(a)
It follows from the above proposition that each non vanishing holomorphic vector field $v$
on $X$ is included in an Abelian Lie algebra $\mf a$ fulfilling~(\ref{descomposicio}).

(b)
Assume that the subalgebra $\mf a$ of $\mf h$ fulfilling the properties of Proposition
(\ref{prop:LieSubA}) has
been fixed and set $\mf a= \langle v_1,\dots,v_s \rangle$. The space of holomorphic 1-forms
on $X$ can be decomposed as
$$
H^0(X,\Omega^1_X)=\langle \alpha^1,\dots, \alpha^s \rangle \oplus \langle \beta^1,\dots,\beta^{p-s} \rangle,
$$
where $p=b_1(X)/2=\dim_{\mbb C} \operatorname{Alb}(X)$, and the forms $\alpha^i$ and 
$\beta^k$ fulfill $\alpha^i(v_j)=\delta^i_j$
and $\beta^k(v_j)=0$ for each $j=1,\dots,s$.

Although the forms $\alpha^i$ are not uniquely determined unless $p=s$, the subspace 
$\langle \beta^1,\dots,\beta^{p-s}\rangle$
is canonically associated to the complex manifold $X$. In particular, it does not depend on 
the choice of the subalgebra $\mf a$.
In fact, $\langle \beta^1,\dots,\beta^{p-s}\rangle$ is just the kernel of the natural morphism 
$$
H^0(X, \Omega^1_X) \longrightarrow H^0(X, \Theta_X)^\vee,
$$
where $\Theta_X$ denotes the sheaf of germs of holomorphic vector fields on $X$.
\end{remarks}

A manifold $X$ will be called {\sl ruled} 
if it is bimeromorphic to a {\sl geometrically ruled} manifold, that is if $X$ is bimeromorphic to the total space 
of a fiber bundle $Y \to B$, which is analytically locally trivial over a complex manifold $B$, has fiber a complex 
projective space $\mathbb C P^k$ and has changes of 
trivialization in the sheaf $\operatorname{PGL}(k+1,\mathcal O_B)$. 
In the case in which $X$ is algebraic this 
is equivalent, by the GAGA theorem, to say that $X$ is bimeromorphic to the trivial bundle 
$\mathbb C P^k\times B$ but this 
equivalence is no longer true in the setting of K\"{a}hler manifolds (cf. \cite{Fuji2}). 
A manifold $X$ will be called {\sl uniruled} if every point $x\in X$ lies in a rational curve $C\subset X$.

\begin{theorem} [Fujiki, \cite{Fuji2}] \label{thm:reglatge}
If a compact manifold $X$ in the class $\mc C$ admits a non-trivial 
holomorphic vector field $v$
with $v_x=0$ for some point $x \in X$, then $X$ is uniruled.
\end{theorem}

\begin{remark}
If the manifold $X$ in Theorem \ref{thm:reglatge} is complex projective then it follows from
the birational classification theorem (Theorem~\ref{thm:classic}) that $X$ is ruled. 
In \cite{Lieb}, Lieberman states that if $X$ is a compact K\"ahler manifold then it is also ruled.
In the particular case when the automorphism group of $X$ includes a $\mbb C^*$-action this
was proved by Carrell and Sommese in \cite{CaSom}. The authors have not been able to
complete Lieberman's argument from \cite{Lieb} in the general compact K\"ahler case. 
\end{remark}

We end this section by recalling the following result of Lieberman concerning the Kodaira dimension,
$\operatorname{kod}(X)$, of $X$ and the dimension of the Lie algebra $\mathfrak{h}_X$.

\begin{theorem}[Lieberman \cite{Lieb}]
Let $X$ be a compact K\"ahler manifold. Then
\begin{equation}\label{thm:lieb2}
\dim_{\mbb C}\mathfrak{h}_X+\mathrm{kod(X)}\leq \dim_{\mbb C}X.
\end{equation}
\end{theorem}

%%%%%%%%%%%%%%%%%      SECTION 2

\section{Locally free actions on K\"{a}hler and Fujiki manifolds}\label{Sec:examples}

Let $X$ be a compact manifold endowed with a locally free $G$-action.
In case $X$ is a K\"{a}hler or Fujiki manifold there are very few possibilities 
for the Lie group $G$.
F. Bosio has remarked the following corollary
of the mentioned theorem of Carrell and Lieberman.

\begin{theorem}[Bosio \cite{Bos}]
Let $G$ be a connected complex Lie group acting locally freely on a compact 
manifold $X$ belonging to the class $\mc C$. Then $G$ is Abelian.
\end{theorem}

\begin{proof}
Theorem~\ref{alpha} implies that the induced action of $G$ on the Albanese torus $\A{X}$
is also locally free.
\end{proof}

Here are some examples of compact complex manifolds naturally endowed with
locally free $\mathbb C^s$-actions.

\begin{example} [Tori]
Complex tori, that is $ \mathbb T^s = \mathbb C^s/\Lambda$ where
$\Lambda\cong\mathbb Z^s$
is a lattice, are the only compact manifolds $X$ with a
locally free $\mathbb C^s$-action such that $s=\dim_{\mbb C} X$.
\end{example}

\begin{example} [Principal torus fibrations]\label{example_fibration}
The total space $X$ of a principal torus bundle
$$
  \begin{array}{ccc}
    \mathbb T^s & \longrightarrow & X \\
     && \downarrow \\
     && B  \\
  \end{array}
$$
where $B$ is a compact complex manifold, is endowed with a natural action of
$\mathbb T^s$ and therefore of $\mathbb C^s$. If $X$ is a K\"{a}hler manifold then
the K\"{a}hler metric can be made invariant by averaging it with respect to the
action by the  compact group $\mathbb T^s$. Hence the base space $B$ is
necessarily a K\"{a}hler manifold. Furthermore, it will be seen below
(cf. Remark \ref{flat_fibration} (a)) that the
principal bundle is flat. Conversely, if $B$ is a
K\"{a}hler manifold and the bundle is flat then the total space is a K\"{a}hler manifold.
This follows from a theorem of Blanchard (cf. \cite{Bla}, Théorème principal II).
\end{example}

\begin{example}[Suspensions over complex tori]\label{example_suspension}
Let $\rho:\Lambda\rightarrow \operatorname{Aut}_{\mbb{C}}(F)$ be a group representation,
where
$\Lambda=\langle \tau_{1},\dots,\tau_{2s}\rangle\cong\mathbb Z^{2s}$
is a lattice of $\mbb C^s$ and $F$ is a compact complex manifold, and denote
$f_{i}=\rho(\tau_i)\in\operatorname{Aut}_{\mbb{C}}(F)$. The {\sl suspension} of the
representation $\rho$ is the manifold $X = F\times_{\Lambda} \mbb C^s$ obtained as the
quotient of the product $F\times\mbb C^s$ by the equivalence relation defined by the
diagonal action of $\Lambda$, that is
$$
(x,z)\sim (f_{i}(x) , z + \tau_{i}), \qquad \text{for}\qquad i = 1,\dots, 2s.
$$
The suspension $X = F\times_{\Lambda} \mbb C^s$ fibers over the torus
$\mbb T^s =\mbb C^s/\Lambda$ with fiber $F$ and the natural $\mathbb C^s$-action on
$F\times\mbb C^s$ commutes with $\Lambda$
inducing a locally free action on $X$. We say that $X$ is a suspension over $\mbb T^s$. 
It is proved in \cite[Theorem 3.19]{Man} that $X$ is a K\"ahler
manifold if and only if $F$ is K\"ahler and there are integers $n_{i}$, for $i = 1,\dots, s$,
such that $f^{n_{i}}_{i}\in\operatorname{Aut}^0_{\mbb{C}}(F)$.
\end{example}

\begin{remark}\label{rem:prob_fujiki}
The above characterizations of K\"{a}hler manifolds
that are principal torus bundles or suspensions over tori stated in the last two examples
rely on a criterium, due to Blanchard \cite{Bla}, for deciding if the total space of a fibration is a K\"{a}hler manifold.
It is not know if Blanchard's criterium can be extended to the class $\mc C$. Therefore we are not able
to formulate such kind of characterizations for Fujiki manifolds. This is the main difficulty to extend 
the structure theorems in Section~\ref{sec:structure} to manifolds in the class $\mc C$.
\end{remark}

From now on $X$ will denote a given compact manifold in the class $\mc C$ and we will use 
the notation introduced in Section \ref{Sec:Prel}. We denote $b_1(X)=2p$. We also assume 
from now on that
$s=\dim_{\mbb C} \mf h/\mf h_{1}$ is strictly positive and we fix a subalgebra
$\mf a\subset \mf h_X$ with $\mf h = \mf h^{1}_X \oplus \mf a$ as well as a basis
$v_{1},\dots,v_{s}$ of $\mf a$. We denote by $\varpi$ the locally free $\mathbb C^s$-action
defined by the choice of the basis and by $\mc F$ the induced foliation.

It follows from 
Theorem~\ref{thm:fuji-lieb} that
there are holomorphic 1-forms $\alpha^{1},\dots,\alpha^{s}$
such that  $\alpha^{i}(v_{i}) = 1$ and $\alpha^{j}(v_{i}) = 0$ if $j\neq i$. The forms
$\alpha^{i}$ generate the cotangent bundle $T^*\mc F$ of $\mc F$ at each point
and we assume that they have also been fixed. The distribution
$\ker\alpha^1\cap\dots\cap\ker\alpha^s$
is integrable and defines a holomorphic foliation
$\mathcal G$ that depends on the choice of  $\alpha^{1},\dots,\alpha^{s}$.
This foliation is transverse and complementary to $\mc F$, and also invariant by the
action of $\mathbb C^s$. Therefore, we have

\begin{proposition}\label{foliation_G}
A foliation $\mc F$ on a compact manifold $X$ in the class $\mc C$ defined by a locally
free $\mbb C^s$-action $\varpi$ admits a holomorphic foliation $\mathcal G$, which is
transverse and complementary to $\mc F$, and invariant by the
action $\varpi$.
\end{proposition}

\begin{remarks} \label{flat_fibration}
(a) This proposition implies in particular that a principal torus bundle whose total space is a
manifold in the class $\mc C$ is necessarily flat (cf. Example \ref{example_fibration}).

(b) Notice that all the leaves of $\mc G$ are biholomorphic. In the case in which the leaves of $\mc G$
are compact they define a fibration of $X$ over a complex torus of dimension $s$ and $X$ 
is a suspension.
\end{remarks}

We denote by $\nu \mc F$ the normal bundle of the foliation $\mc F$, i.e.
$\nu \mc F = TX/T\mc F$.
The existence of the $\mbb C^s$-invariant foliation $\mc G$ implies that, 
at each point of $X$, we can find local coordinates
$(t^1, \dots , t^s, z^1, \dots , z^{n-s})$ with the properties
\begin{equation}\label{local_coordinates}
v_{i} = \frac{\partial}{\partial t^{i}} \qquad \text{and} \qquad \alpha^{i} = d t^{i}.
\end{equation}
In these coordinates $\mc  F$ is defined by $z^{k} = const.$ and
$\mc  G$ by $t^{i} = const$.
As a consequence, we obtain:

\begin{proposition}\label{prop:canonical_bundle}
There is a natural isomorphism between the canonical bundle $K_{X}$ of $X$
and the determinant $\det \nu^* \mc F$ of the conormal bundle 
$\nu^* \mc F$ of $\mc F$.
\end{proposition}

\begin{proof}
If $(t_{i}^1, \dots , t_{i}^s, z_{i}^1, \dots , z_{i}^{n-s})$ are local coordinates 
fulfilling (\ref{local_coordinates}) and we define
$\eta_{i} = \alpha^1\wedge\dots\wedge\alpha^s\wedge dz_{i}^1\wedge\dots\wedge dz_{i}^{n-s}$,
then
$\eta_{j} = g_{ji} \eta_{i}$ with
$$
g_{ji} = \det \bigg(\frac{\partial z_{j}^k}{\partial z_{i}^l} \bigg).
$$
Hence the cocycle $(g_{ji})$ defines the line bundle $K_{X}$ as well as
$\det \nu^* \mc F$.
\end{proof}

We then obtain:

\begin{corollary}\label{cor:canonical_bundle}
Assume that $X$ is a suspension $F\times_\Lambda \mbb C^s$ over a torus $T = \mbb C^s/ \Lambda$.
Then one has $\operatorname{kod}(X)=\operatorname{kod}(F)$ and $c_1(F) = \iota^* c_1(X)$
where $\iota\!: F\hookrightarrow X$ is the identification of $F$ with a fiber of the projection $X\to T$.
\end{corollary}

\begin{proof}
Notice that the tangent bundle $TF$ of $F$ is canonically isomorphic to the restriction of $\nu\mc F$ to $F$.
As the Kodaira dimension is invariant by finite coverings we can assume, using \cite[Theorem 3.19]{Man}, 
that the monodromy of the suspension has values in $\operatorname{Aut}^0_{\mbb C}(F)$.
By Proposition~\ref{prop:canonical_bundle} there are isomorphisms 
$K_X^{\otimes m}\cong (\det \nu^* \mc F)^{\otimes m}$ for every $m\in \mbb Z$. Consequently 
the spaces $H^0(X, K_X^{\otimes m})$ are naturally identified with the subspaces
of $H^0(F, (\det \nu^* \mc F)^{\otimes m}) = H^0(F, K_F^{\otimes m})$ invariant by the monodromy.
Ueno proved \cite[Corollary 14.8]{Ue}  that $\operatorname{Aut}^0_{\mbb C}(F)$ acts 
trivially on the pluricanonical section spaces $H^0(F, K_F^{\otimes m})$. 
Therefore $h^0(X, K_X^{\otimes m})= h^0(F, K_F^{\otimes m})$ for every $m$. 

The relation between the first Chern classes follows from the adjunction formula.
\end{proof}

Let  $\Omega^k_X$ be the sheaf of germs of holomorphic k-forms on $X$.
As remarked above $p={b_{1}(X)}/{2}\geq s$.
We can choose holomorphic $1$-forms $\beta^{1},...,\beta^{p-s}$ on $X$
such that
\begin{equation}\label{eq:base_formes}
\{ \alpha^{1},...,\alpha^{s},\beta^{1},...,\beta^{p-s} \}
\end{equation}
is a basis of $H^0(X,\Omega^1_X)$ and that $\beta^j(v_i)=0$ for $i=1,\dots,s$ and
$j=1,\dots,p-s$.
Since $\beta^{j}$ are closed they are in fact {\sl basic} forms with
respect to $\mathcal F$. We recall that a differential
form $\gamma$ is said to be basic with respect to a foliation if it fulfills $i_{w}\gamma =0$ and
$i_{w}d\gamma = 0$ for each local vector field $w$ tangent to the foliation.
We shall denote by $\Omega^k_{X/\mc F}$ the subsheaf of $\Omega^k_X$ of
those k-forms that are basic with respect to $\mc F$.
The following decomposition of global holomorphic forms on $X$ will be useful
in the sequel.

\begin{lemma}\label{decomposition}
Each global holomorphic $k$-form $\gamma$ on $X$ can be written in a unique way
as a sum
$$
\gamma=\sum_{0\leq |I|\leq \min(k,s)} \alpha^{I} \wedge \gamma^I
$$
with $I = (i_{1},\dots,i_{m})$, $1 \leq i_{1} < \dots < i_{m}\leq s$ and $|I| = m$,
and where $\alpha^{I}=\alpha^{i_1} \wedge \dots \wedge \alpha^{i_m}$ and
$\gamma^I$ is a holomorphic $(k-m)$-form \emph{basic} with respect
to $\mathcal{F}$. In particular we have
$$
H^0(X,\Omega^k_X)\cong \bigoplus_{i=0}^{\min(k,s)} \bigwedge^{i} \langle \alpha^{1},...,\alpha^{s} \rangle_{\mbb C} \otimes H^0(X,\Omega^{k-i}_{X/\mc F}).
$$
\end{lemma}

\begin{proof}
Let $\gamma$ be a global holomorphic $k$-form on $X$ and set
$\ell = \min(k,s)$. For each $I = (i_{1},\dots,i_{m})$, with $m\leq k$ and
$1 \leq i_{1} < \dots < i_{m}\leq s$, and a given $k$-form $\omega$ we denote
$\omega^{I} = i_{v_{i_{m}}}\circ \dots \circ i_{v_{i_{1}}}\omega$. We put
$\gamma_{0} =\gamma$ and we define $\gamma_{1},\dots,\gamma_{\ell}$
recursively by
$$
\gamma_{m+1}=\gamma_{m} - \sum_{|I| = \ell -m} \alpha^{I} \wedge \gamma_{m}^I.
$$
Notice here that the $(k-\ell + m)$-forms $\gamma_{m}^{I}$ are holomorphic global forms
and thus they are closed.
Moreover, by construction they fulfill $i_{v}\gamma_{m}^{I}=0$ for each vector $v$ tangent to
$\mc F$. Therefore they are basic with respect to $\mc F$ and
$$
\gamma = \sum_{m=1}^\ell \sum_{|I|=m} \alpha^{I}\wedge \gamma_{\ell-m}^{I} + \gamma_{\ell}
$$
is the required decomposition.
\end{proof}

%%%%%%%%%%%%%%%%% SECTION 3

\section{Projective manifolds with nonvanishing tangent fields} \label{s:proj}

In this section we will assume that $X$ is a complex projective manifold. 
Our study in this algebraic context becomes a continuation of the classical study of complex
projective manifolds with holomorphic tangent vector fields, which led to their birational
classification summarized in Theorem \ref{thm:classic}. The Albanese mapping and 
Poincar\'e's reducibility theorem lead to a biholomorphic classification of these manifolds 
which is simpler than in the K\"{a}hler case. 

Following Theorem~\ref{thm:fuji-lieb} and Remark~\ref{rem:camps} about the structure of the Lie 
algebra of holomorphic tangent vector fields, the existence of nonvanishing tangent 
vector fields is equivalent to the existence of a holomorphic locally free $\mbb C^r$--action
on the manifold. We will work directly with these locally free actions. The basic fact 
beyond the classical theory is:

\begin{proposition} \label{prop:susp_proj}
Let $X$ be a complex projective manifold, admitting a locally free holomorphic
$\mbb C^r$-action. Then $X$ is a suspension $F\times_{\Lambda} \mbb C^s$ over 
an Abelian variety of dimension $s\geq r$, with fiber a connected projective manifold
$F$.
\end{proposition}

\begin{proof}
As $X$ is closed K\"ahler, the $\mbb C^r$-action can be extended to a locally free holomorphic
action of maximal rank $s=\dim_{\mbb C} T_X=\dim \mathfrak h_{X}/\mathfrak h_{X}^1$ as
we have recalled in Section \ref{Sec:Prel}. 

The Albanese torus $\A{X}$ of $X$ is an Abelian variety, therefore  the subtorus 
$T_X \subset \text{Alb}\,(X)$ in the 
Fujiki-Lieberman structure Theorem~\ref{thm:fuji-lieb}
is also an Abelian
variety, and by Poincar\'e's Reducibility Theorem (see \cite[\S 5.3]{BL}) it has a complementary
Abelian subvariety $Z$ such that $T_X \cap Z$ is finite, $T_X+Z=\A{X}$, and the addition law induces
an isogeny $c \!: \A{X} \rightarrow T_X \times Z$.

The composition of the Albanese morphism of $X$ with this isogeny and the natural projection
$$
X \stackrel{\phi}{\longrightarrow} \A{X} \stackrel{c}{\longrightarrow} T_X \times Z \stackrel{p_1}{\longrightarrow} T_X
$$
maps the tangent subspace of $X$ generated by the $\mbb C^s$-action
isomorphically into the tangent space of $T_X$. Consequently $X$ is a suspension over $T_X$,
with parallel transport given by the $\mbb C^s$-action and possibly disconnected
fibers. The Stein factorization of $X \to T_X$ is thus unramified, and yields a 
suspension $X \to N_X$ with connected fiber $F$, over an Abelian variety $N_X$  
isogenous to $T_X$.
\end{proof}

If $X$ is not ruled, it follows from Theorems \ref{thm:classic} and \ref{thm:fuji-lieb}
that $\operatorname{Aut}^0_{\mbb C}(X)$ is an Abelian variety of dimension $s =\dim \mathfrak 
h_{X}$, isogenous to $N_X$. The suspension morphism $X \to N_X$ shows that 
the natural action of $\operatorname{Aut}^0_{\mbb C}(X)$ on $X$ 
is holomorphically injective, and the above quoted theorems
make Proposition~\ref{prop:susp_proj} essentially equivalent
to Carrell's classification theorem for such actions of Abelian varieties
(cf. \cite[ Thm. 8]{Car71}). 

The classification of projective manifolds with locally free holomorphic $\mbb C^s$-actions
follows from Proposition~\ref{prop:susp_proj} and Theorem~\ref{thm:reglatge}. Putting them together yields:

\begin{theorem} \label{thm:susp_proj}
Let $X$ be a complex projective manifold admitting a locally free holomorphic 
$\mbb C^r$-action.
\begin{enumerate}[\rm (a)]
\item If $X$ is not uniruled, it is a quotient $(F \times T)/\Gamma$, with $T$
an Abelian variety of dimension $s\ge r$, $F$ a projective manifold with
no holomorphic tangent vector fields and $\operatorname{kod}(F) = \operatorname{kod}(X)$,
and $\Gamma$ a finite Abelian group of biholomorphisms operating freely
on $F \times T$.
\item If $X$ is uniruled, it admits a finite, Abelian, \'etale cover $X'$ which is 
a suspension over an Abelian variety $T$ of dimension $s\ge  r$, with uniruled fiber 
$F$ such that any holomorphic tangent vector field on $F$ has zeros. The monodromy of 
the suspension $X' \to T$ has values in the group $\operatorname{Aut}^0_{\mbb C}(F)$.
\end{enumerate}
\end{theorem}

\begin{remark} This Theorem proves the algebraic part of  Theorem~\ref{teo:mainmain}.
Its proof mirrors that of Proposition~\ref{prop:primera_reduccio} which 
deals of the general case of K\"{a}hler manifolds. As examples in Section~\ref{sec:structure} 
show, the fact that 
Poincar\'e's reducibility theorem is no longer valid in the K\"{a}hler setting forces the introduction 
of deformations of the complex structure. 
\end{remark}

\begin{proof} 
We may assume that the locally free action on $X$ is of maximal rank $s$.
By Proposition~\ref{prop:susp_proj}, $X$ is a suspension over an $s$--dimensional 
Abelian variety $N_X$, isogenous to $T_X$, and with fiber a connected projective 
manifold $F_1$.

Let the suspension structure on $X$ be given by $N_X = \mbb C^s/ \Lambda$,
with $\Lambda = \langle \tau_1, \dots , \tau_{2s} \rangle$ a cocompact lattice defining
the torus $N_X$, $\rho : \Lambda \rightarrow \operatorname{Aut}_{\mbb C}(F)$ the monodromy
of the suspension, and $f_i = \rho(\tau_i)$ the generators of the monodromy.

By \ref{example_suspension}, the fact that the total space $X$ of the suspension is K\"ahler
implies that all the monodromy automorphisms have a finite power $f^{n_i}_i \in  
\operatorname{Aut}^0_{\mbb C}(F)$,
i.e. there is a finite \'etale Abelian covering $T_1 \rightarrow N_X$ such that the
pullback of $X$ over $T_1$ is a finite \'etale Abelian covering $X_1 \rightarrow X$,
and also a suspension over $T_1$ with fiber $F_1$ and monodromy automorphisms in
$ \operatorname{Aut}^0_{\mbb C}(F_1)$.

If $F_1$ has no nonvanishing holomorphic tangent vector fields it can be checked
that we have reached the sought cover: 
\begin{enumerate}[\rm (i)]
\item $\operatorname{kod}(F_1) = \operatorname{kod}(X)$ by Corollary~\ref{cor:canonical_bundle}.
\item If $X$ is not uniruled neither is $X_1$ nor $F_1$. Indeed, as all the fibers of $X_1\to T_1$
are isomorphic, if $F_1$ were uniruled there would be a rational curve passing by each given point of $X_1$.
Now Fujiki's Theorem~\ref{thm:reglatge} implies that $\mf h^1_{F_1} =0$. 
Therefore the suspension monodromy
has values in $\operatorname{Aut}^0_{\mbb C}(F_1)= \{ \text{Id} \}$,
i.e. $X_1$ is biholomorphic to $F_1 \times T_1$.
\item If $X$ is uniruled so is its \'etale cover $X_1$. Since rational
curves map to points in the suspension basis $T_1$, the fiber $F_1$ is uniruled
as well. 
\end{enumerate}

If $F_1$ has nonvanishing holomorphic tangent vector fields, we apply an iterative
process: let $\mf h_{F_1}$
be the Lie algebra of holomorphic vector fields in $F_1$, and apply Proposition
\ref{prop:touzet} to choose an Abelian subalgebra $\mf a_{F_1}$ corresponding 
to a $\mbb C^r$--action of maximal rank in $F_1$, such that $\mf a_{F_1}$ 
is central in $\mf h_{F_1}$.

As the monodromy of the suspension $X_1 \to T_1$ lies in $\operatorname{Aut}^0_{\mbb C}(F_1)$
it must commute with the central algebra $\mf a_{F_1}$. Therefore its 
tangent vector fields may be extended by parallel transport from the fiber $F_1$ 
to the suspension total space $X_1$ as nonvanishing tangent 
vector fields. These extensions, to be denoted $\langle v_1, \dots, v_r \rangle$,
commute with the horizontal vector fields
$\langle h_1, \dots, h_s \rangle$ which are the lifts to $X_1$ of $\mf h_{T_1} = H^0(T_1, \Theta_{T_1})$
defining the suspension's parallel transport. Thus $\mf a_{X_1}= \langle v_1, \dots , v_r, h_1\dots
h_s \rangle$ is an Abelian subalgebra of $\mf h_{X_1}$, spanned by $r+s$ linearly independent
nonvanishing vector fields.
It has maximal rank among such algebras: the algebra $\mf h_{X_1} $ is the direct sum 
$$\mf h_{X_1}  = \langle h_1,\dots, h_s \rangle \oplus H^0(X_1, \mc {V}ert)$$ 
where $\mc {V}ert$ is the sheaf spanned  by the vector fields vertical with respect to the projection
$X_1\to T_1$. As the horizontal fields $h_1,\dots, h_s$ are central in $\mf h_{X_1}$,
the morphism  $H^0(X_1, \mc {V}ert)\to \mf h_{F_1}$ given by restriction to a fiber
is injective. This induces a natural inclusion of Lie algebras
$$
\mf h_{X_1}  = \langle h_1,\dots, h_s \rangle \oplus H^0(X_1, \mc {V}ert) \subset
\langle h_1,\dots, h_s \rangle \oplus \mf h_{F_1}.
$$
If there existed an Abelian Lie subalgebra $\tilde{\mf a}\subset \mf h_{X_1}$ formed 
by nonvanishing tangent fields, with $\dim\tilde{\mf a}> r+s$, then its inclusion
in $\langle h_1,\dots, h_s \rangle \oplus \mf h_{F_1}$ would have an intersection with $\mf h_{F_1}$
of dimension $>r$. This contradicts the maximality of $\mf a_{F_1}.$

By Proposition~\ref{prop:susp_proj}, the algebra $\mf a_{X_1}$ determines a suspension $X_1 \to N_1$,
over an $(r+s)$--dimensional Abelian variety, with a connected fiber 
$F_2$. 
Again by \ref{example_suspension}, the monodromy automorphisms of this suspension have a
finite power in $\operatorname{Aut}^0_{\mbb C}(F_2)$. Consequently, there is a finite
\'etale Abelian covering $T_2 \rightarrow N_1$ such that the
pullback $X_2$ of $X_1$ over $T_2$ is a suspension with fiber $F_2$ and
monodromy in $\operatorname{Aut}^0_{\mbb C}(F_2)$. 

The induced map $X_2 \to X_1$ is also a finite \'etale Abelian cover.
Moreover, the composition of coverings $X_2 \to X_1 \to X$ is still
regular and Abelian. The reason for its regularity is that the automorphisms
of the cover $X_1\to X$ are realized by integration up to time 1 in $X_1$ of 
suitable linear combinations of the horizontal fields $h_1, \dots, h_s$.
These tangent fields have become parallel transport vector fields
of the suspension $X_1 \to T_1$ by our choice of algebra $\mf h_{X_1}$,
and the isogenous base change $T_2 \to T_1$ lifts canonically to 
$X_2$ the vector fields $h_1, \dots, h_s$, thus also the integration
up to time 1 of their linear combinations. The Abelianity of the 
cover follows from the fact that the vector fields $v_1, \dots , v_r, h_1\dots h_s$
commute in $X_1$, hence 
also in $X_2$, and the automorphisms of both covers $X_1\to X$, $X_2\to X_1$ 
are defined by integration up to time 1 of suitable linear combinations
of them.

The iteration step concludes here. We have obtained a suspension 
$X_2 \to T_2$, with $X_2$ a finite, Abelian \'etale cover of $X$,
fiber $F_2$ with $\operatorname{kod}(F_2) = \operatorname{kod}(X_2)= \operatorname{kod}(X)$,
and monodromy in $\operatorname{Aut}^0_{\mbb C}(F_2)$. 

Repetition of the above procedure yields a sequence
of suspensions $X_i \to T_i$ with fiber $F_i$, again with $X_i$ Abelian \'etale
covers of $X$, $\operatorname{kod}(F_i) = \operatorname{kod}(X)$ and 
monodromy in $\operatorname{Aut}^0_{\mbb C}(F_i)$. As the dimension of the 
basis tori $T_2, T_3, \dots$ grows strictly the iterative process must
reach a final suspension $X_f \to T_f$ such that its fiber $F_f$ does not
have any nonvanishing vector field, i.e. no holomorphic locally free
$\mbb C^r$--action. 

If $X$ is not uniruled neither is $X_f$ nor $F_f$. By 
Fujiki's Theorem \ref{thm:reglatge} we have $\mf h^1_{F_f}=0$. Hence
the fiber $F_f$ cannot have any
tangent vector field so $\operatorname{Aut}^0_{\mbb C}(F_f)= \{ \text{Id} \}$
and the trivial suspension $X_f=F_f \times T_f$ has been reached. 

If $X$ is uniruled, so are all its \'etale covers up to $X_f$, and 
the fiber $F_f$ as well because rational curves in $X_f$ must map to a point
in $T_f$. 
\end{proof}

We can use the structure Theorem~ \ref{thm:susp_proj} to classify algebraic manifolds 
with sufficiently many
vector fields and to establish for them Ueno's conjecture~\ref{conj:kollar}. This is done
in Corollaries~\ref{c:clasif_proj} and \ref{c:k=0} below as specialization of the arguments in the 
K\"{a}hler case.

%%%%%%%%%%%%%%%%%%%%%     SECTION 4

\section{Tangential deformations of locally free holomorphic actions}\label{Sec:TanDef}

Through this section $X$ will be an arbitrary compact complex manifold,
i.e. not necessarily K\"{a}hler or belonging to the Fujiki class.
Let us assume that $X$ is endowed with a holomorphic action
$\varpi: G\times X\rightarrow X$, where $G$ is a connected complex Lie group.
Such an action is defined by a representation
\begin{equation}\label{action}
\rho: G\longrightarrow \operatorname{Aut}_{\mbb{C}}(X)
\end{equation}
from $G$ into the group $\operatorname{Aut}_{\mbb{C}}(X)$ of holomorphic
automorphisms of $X$. When $\varpi$ is locally free, i.e. it has discrete isotropy groups, 
then it induces a holomorphic foliation $\mathcal F$ on $X$ whose leaves are the
orbits of the action.

Let  $\mathcal F^{tr}$ denote the transversely holomorphic foliation obtained from
$\mathcal F$ by forgetting the complex structure along the leaves. One can consider
deformations of the holomorphic foliation
$\mathcal F$ that keep fixed its transversal type; that is, holomorphic deformations
$\mathcal F_{r}$ of $\mathcal F$ such that the transversely
holomorphic foliation $\mathcal F^{tr}_{r}$ coincides with
$\mathcal F^{tr}$. This type of deformations were studied in \cite{G-M} and \cite{GN},
where they are called $f$-deformations. In \cite{G-M}, Gómez-Mont studies the
space of infinitesimal $f$-deformations of holomorphic foliations, which is naturally
identified to $H^1(X, T{\mathcal F})$, i.e. the first cohomology
group of $X$ with values in the tangent bundle $T{\mathcal F}$
of $\mathcal F$. The existence of a {\sl versal} or Kuranishi space for $f$-deformations,
i.e. a germ of analytic space $(\K^{f},0)$ parametrizing a versal family of $f$-deformations
of $\mathcal F$, is proved in \cite{GN}. The tangent space of $(\K^{f},0)$ at $0$
is $H^1(X, T{\mathcal F})$.

Families of $f$-deformations of $\mathcal F$ can be viewed as unfoldings of $\mathcal F$,
that is families of complex structures on $X$ for which $\mathcal F^{tr}$ becomes holomorphic.
More precisely, a family of $f$-deformations of $\mathcal F$ parametrized by a germ of
analytic space $(R, r_{0})$ is given by a proper and flat morphism
$\pi:X_{R}\rightarrow R$, an identification $\iota: X\hookrightarrow X_{R}$ of $X$
with the fiber $X_{r_{0}}=\pi^{-1}(r_{0})$ and
a holomorphic foliation $\mathcal F_{R}$ on $X_{R}$ of the same codimension than $\mathcal F$,
which is transverse to the projection $\pi$ and such that the restriction
$\mathcal F_{r_{0}}$ of $\mathcal F_{R}$ to $X\equiv X_{r_{0}}$ coincides
with $\mathcal F$. Such a family of $f$-deformations will be denoted by
$(X_{R}, \iota, \pi, \mathcal F_{R})$.

Notice however that $f$-deformations $\mathcal F_{r}$ of $\mathcal F$ need not
to be equivariant, in the sense that $\mathcal F_{r}$ are not necessarily defined
by an action. By this reason we give the following definition, which is inspired in
the notion of equivariant deformations of complex manifolds introduced
by Cathelineau in \cite{Cat}.

\begin{definition}
Let $X$ be a compact complex manifold endowed with a locally free holomorphic action
$\varpi: G\times X\rightarrow X$ and let $\mathcal F$ be the foliation defined by $\varpi$.
By a {\sl family of tangential deformations}
of $\varpi$ we mean
a family of $f$-deformations $(X_{R}, \iota, \pi, \mathcal F_{R})$ of $\mathcal F$
and a holomorphic action $\varpi_{R}$ of $G$ on $X_{R}$ fulfilling the following properties:
\begin{enumerate}[(i)]
\item $\varpi_{R}$ extends the action $\varpi$ of $G$ on $X\equiv \iota(X) = X_{r_{0}}$,
\item $\varpi_{R}$ preserves the fibers of $\pi$ and
\item the orbits of the action are tangent to $\mathcal F_{R}$.
\end{enumerate}
Then $\varpi_{R}$ induces a locally free $G$-action $\varpi_{r}$ on $X_{r}$ for each $r\in R$
and the leaves of $\mathcal F_{r} = \mathcal F|_{X_{r}}$ are the orbits of $\varpi_{r}$.

Such a family will be denoted by $(X_{R}, \iota, \pi, \varpi_{R})$ and each pair $(X_r, \varpi_{r})$,
with $r\in R$, will be called a {\sl tangential deformation} of $(X_0, \varpi_{0})\cong (X, \varpi)$.
We will just write $(X_{R}, \varpi_{R})$  instead of  $(X_{R}, \iota, \pi, \varpi_{R})$ when
there is no danger of confusion.
\end{definition}

\begin{definition}
Two families $(X_{R}, \iota, \pi, \varpi_{R})$ and $(X'_{R}, \iota', \pi', \varpi'_{R})$ of tangential 
deformations of $(X, \varpi)$,  parameterized by the same germ of analytic space $R$, are 
said to be isomorphic if there is a $G$-equivariant
biholomorphism $\phi: X_R\rightarrow X'_R$ fulfilling $\phi\circ\iota = \iota'$ and $\pi'\circ\phi = \pi$.
\end{definition}

A family of tangential deformations $(X_{R}, \varpi_{R})$ of $(X, \varpi)$ is called
{\sl versal} if it fulfills the following property: for any other family $(X_{S}, \varpi_{S})$
of tangential deformations of $(X, \varpi)$ there is an analytic morphism of germs
of analytic spaces
$\varphi: S \rightarrow R$ such that
\begin{enumerate}[(i)]
\item the pull-back family $(\varphi^*(X_{R}), \varphi^*(\varpi_{R}))$ parameterized by $S$ is 
isomorphic to $(X_{S}, \varpi_{S})$, and
\item the tangent map $d_{s_0}\varphi$ of $\varphi$ at the distinguished point $s_0$ of $S$ is unique.
\end{enumerate}
If such a versal family exists then it is unique up to isomorphism.
\medskip

Notice now that the sheaf of germs of holomorphic vector fields tangent to $\mathcal F$
is a $G$-sheaf and therefore we can consider the equivariant cohomology of $X$
with values in the tangent bundle $T\mathcal F$ of $\mc F$, that we denote by $H_{G}^*(X,T\mathcal F)$.
Let $\Omega^{0,q}(X, T\mathcal F)$ be the Frechet space
of $(0,q)$-forms on $X$ with values in $T\mathcal F$ and denote by
$C_h^p(G, \Omega^{0,q}(X, T\mathcal F))$ the space of holomorphic $p$-cochains
with values in the $G$-module $\Omega^{0,q}(X, T\mathcal F)$.
Then $H_{G}^*(X,T\mathcal F)$ is the cohomology of the double complex
\begin{equation}\label{bicomplex}
((C_h^p(G, \Omega^{0,q}(X, T\mathcal F)))_{p,q}, \delta, \overline\partial)
\end{equation}
where $\delta$ denotes the Eilenberg-Maclane coboundary operator and
$\overline\partial$ is the usual delta-bar operator. As usual, infinitesimal
tangential deformations of $(X, \varpi)$ can be defined as isomorphism classes of 
families parameterized by the double
point $D=\{\ast, \mbb C(t)/(t^2)\}$. The space of infinitesimal deformations is naturally
identified to $H_{G}^1(X,T\mathcal F)$. Moreover, to each family of tangential deformations
$(X_{R}, \varpi_{R})$ of $(X, \varpi)$ there is associated in a natural way a linear map
$$
\kappa: T_{r_0}R \longrightarrow H_{G}^1(X,T\mathcal F)
$$
which is called the Kodaira-Spencer map of the family  (cf. \cite{Cat}).

The following theorem states the existence of a versal family and it can be proved 
in the same lines as
the corresponding statement for equivariant deformations of complex manifolds
(cf. \cite[Th\'{e}or\`{e}me 1 and Proposition 1]{Cat}).

\begin{theorem}\label{teo:kuranishi_equivariant}
Let $X$ be a compact complex manifold endowed with a locally free holomorphic action
$\varpi: G\times X\rightarrow X$ and let $\mathcal F$ be the induced
holomorphic foliation. Then  there
is a germ of analytic space $(\K_{G}^t, 0)$ parametrizing a versal family of
tangential deformations of the action $\varpi$. The tangent space to $\K_{G}^t$ at $0$ is
isomorphic to $H_{G}^1(X,T\mathcal F)$ and it is
naturally identified to the space of infinitesimal
tangential deformations of $\varpi$.
\end{theorem}

The following proposition is a general fact in deformation theory and can be deduced easily
from the existence of the versal space of tangential deformations.

\begin{proposition}\label{prop:criteri_versal}
Let $(X_{R}, \varpi_{R})$ be a family of tangential deformations of $(X, \varpi)$ whose parameter
space $R$ is smooth and such that the corresponding Kodaira-Spencer map $\kappa$ is
an isomorphism. Then the family of deformations $(X_{R}, \varpi_{R})$ is versal.
\end{proposition}

Associated to the double complex (\ref{bicomplex}) there is an spectral sequence converging to
$H_{G}^*(X,T\mathcal F)$ and whose second term is
$E_{2}^{p,q} = H_{h}^p(G, H^q(X,T\mathcal F))$. Here the index $h$ denotes the cohomology
of holomorphic cochains. In particular one has the following exact sequence that is useful in
the computation of $H_{G}^1(X,T\mathcal F)$
\begin{equation}\label{eq:ShortSequence}
\begin{split}
0\rightarrow H^1_{h}(G, H^0(X,T\mathcal F)) & \rightarrow H_{G}^1(X,T\mathcal F)
\stackrel{\chi}{\rightarrow} H^1(X,T\mathcal F)^G \rightarrow \\
& \rightarrow H^2_{h}(G, H^0(X,T\mathcal F))  \rightarrow H_{G}^2 (X,T\mathcal F).
\end{split}
\end{equation}

%%%%%%%%%%%%%%%%%%%%%% SECTION 5

\section{Tangential deformations of locally free $\mathbb C^s$-actions on K\"{a}hler 
and Fujiki manifolds}\label{sec:versal_family}

The aim of this section is to construct the versal family of tangential deformations of a
locally free $\mathbb C^s$-action $\varpi$ on a given compact K\"{a}hler or Fujiki manifold $X$
and to study the properties of that family.

With this purpose we consider the exact sequence (\ref{eq:ShortSequence})
when $G$ is the Abelian group $\mbb C^s$. In that case $T\mathcal F$ is trivial as a $G$-bundle.
So $G= \mbb C^s$ acts trivially on $H^0(X, T\mathcal F) \equiv \mbb C^s$ and
$H^1_h(G, H^0(X, T\mathcal F))$ is just the space
of holomorphic group homomorphisms from $G$ into $\mbb C^s$, that is the space
$\operatorname{End}\mathbb C^s$ of $\mbb C$-linear endomorphisms of $\mbb C^s$.
Moreover, $G= \mbb C^s$ acts also trivially on $H^1(X, T\mathcal F) =
H^1(X, \mathcal O_X)^s\cong H^0(X,\Omega^1_X)^s$, where the last isomorphism 
is $\mbb C$-antilinear.

Therefore, the first terms of sequence (\ref{eq:ShortSequence}) can be written
\begin{equation}\label{eq:AbelianShortSequence}
0\rightarrow \operatorname{End}\mathbb C^s \stackrel{\tau}\rightarrow H_{G}^1(X,T\mathcal F)
\stackrel{\chi}{\rightarrow} H^0(X,\Omega^1_X)^s \rightarrow \dots
\end{equation}
We will see below that the morphism $\chi$ is surjective and therefore
that $H_{G}^1(X,T\mathcal F)$ is isomorphic to the direct sum
$\operatorname{End}\mathbb C^s \oplus H^0(X,\Omega^1_X)^s$.  This motivates the
construction of the family of tangential deformations of $(X, \varpi)$ that follows.
\bigskip

As above $v_1,\dots, v_s$ will denote the fundamental vector fields of the action $\varpi$
corresponding to the canonical basis of $\mbb C^s$, $\alpha^1, \dots,\alpha^s$ will 
be holomorphic 1-forms fulfilling $\alpha^i(v_j) = \delta^i_j$ and 
$\{ \alpha^{1},...,\alpha^{s},\beta^{1},...,\beta^{p-s} \}$ a basis of $H^0(X,\Omega^1_X)$
where $\beta^k$ are basic forms. Although the vector fields $v_i$ are naturally associated 
to the action $\varpi$, the 1-forms $\alpha^i$ are not uniquely determined. If $\tilde\alpha^i$ 
is another choice then
$\tilde\alpha^i-\alpha^i$ are basic forms, hence vanishing on $T\mathcal F$. 

We denote $\Xi = \operatorname{End}\mathbb C^s \oplus H^0(X,\Omega^1_X)^s$.
Given an element $r=(C,\theta)\in \Xi$,
where  $\theta = (\theta^{1}, \dots,\theta^{s})$ and
$\theta^{i} = \sum_{j=1}^s a_{j}^{i} \alpha^{j} + \sum_{k=1}^{p-s} b_{k}^{i}\beta^{k}$
with $a_{j}^{i}, b_{k}^{i}\in \mbb C$,
we set
\begin{equation}\label{new-alpha}
\alpha_{r}^{i}= \alpha^{i} + \overline\theta^{i} = \alpha^{i} +
\sum_{j=1}^s \overline a_{j}^{i} \overline\alpha^{j} + \sum_{k=1}^{p-s} \overline b_{k}^{i}\overline\beta^{k}
\quad \mbox{for} \quad i=1,\dots,s.
\end{equation}
We denote $A = (a_{j}^{i})$ and we define the matrix $M_{A}$ by
$$
M_{A} =
\left(
\begin{array}{cc}
I & -A  \\
-\overline A & I   \\
\end{array}
\right)
$$

\begin{definition}\label{def:R}
Let $R_0$ be the open subset of  $H^0(X,\Omega^1_X)^s$ of those elements 
$\theta = (\theta^{1}, \dots,\theta^{s})$ for which $\det M_{A} \neq 0$ and denote
by $R$ the open subset of $\Xi$ whose elements are the pairs
$r= (C,\theta)$ fulfilling
\begin{enumerate}[(i)]
\item $\det C \neq 0$ and
\item $\det M_{A} \neq 0$,
\end{enumerate}
that is, $R=\operatorname{GL}(s,\mathbb C)\times R_0$. The set $R$ is the complementary 
of an affine real algebraic variety and its tangent space at $(\mathrm{id},0)$ is naturally identified 
to $\Xi = \operatorname{End}\mathbb C^s \oplus H^0(X,\Omega^1_X)^s$.
\end{definition}

The space $R$ will be the parameter space of the versal family of tangential deformations 
of $(X,\varpi)$. Moreover, in the case in which $X$ is a K\"{a}hler manifold, that family will be 
versal at each point of $R$.

Condition (ii) in the above definition implies that $\alpha_{r}^{1},\dots,\alpha_{r}^{s}$ are 
linearly independent 
at each point. Hence the forms $\alpha_{r}^{i}$
generate a real subbundle $Q_r$ of $T^*X^{\mathbb C}$.
One has $Q_r\cap\nu^*\mc F =0$, so there are well-defined complex valued smooth
vector fields $w_{1},\dots w_{s}$ that are tangent to $\mc F$ and are determined by the conditions
\begin{equation}\label{new-vectors}
\alpha_{r}^{i}(w_{j}) = \delta_{j}^{i} \quad \mbox{and} \quad \alpha_{r}^{i}(\overline{w}_{j}) = 0.
\end{equation}
We set 
\begin{equation}
\tilde w_i = v_i-\sum_{j=1}^s a^j_i \overline{v}_j \quad \mbox{for} \quad i=1,\dots,s,
\end{equation}
and we notice that
\begin{equation}\label{express-vectors}
\langle w_{1},\dots, w_{s}\rangle_{\mbb C}=
\big\langle \tilde w_{1},\dots, \tilde w_{s} \big\rangle_{\mbb C}\cong \mbb C^s.
\end{equation}
Then we define 
\begin{equation}\label{new-vectors-last}
v^r_i = Cw_i. 
\end{equation}

In this situation one has

\begin{proposition}\label{prop:integrabilitat}
Let an element $r=(C,\theta)$ of $R$ be given and set
$N_r^{1,0} = Q_r \oplus \nu^*\mc F^{1,0}$. Then 
$T^*X^{\mathbb C} = N_r^{1,0}\oplus\overline{N_r^{1,0}}$ and the almost complex structure
so defined is integrable defining a complex structure $X_{r}$. Moreover, the vector fields
$v_{r}^{i}$ are holomorphic and induce a locally free
$\mbb C^{s}$-action $\varpi_{r}$ on $X_{r}$, which is a tangential deformation of $\varpi$.
This construction defines a family $(X_R,\varpi_R)$ of tangential deformations of $(X,\varpi)$
parametrized by $R$.
\end{proposition}

\begin{proof} Condition $T^*X^{\mathbb C} = N_r^{1,0}\oplus\overline{N_r^{1,0}}$ is
equivalent to $\det M_{A} \neq 0$. Since $\alpha^i_r= \alpha^i + \overline\theta^i$ are closed, 
the almost complex structure on $X$ defined by $N_r^{1,0}$ is integrable and 
induces a complex structure $X_r$ by Newlander-Nirenberg's theorem.

Since $v_1,\dots,v_s$ are holomorphic on $X$ we see, using~(\ref{express-vectors}),
that the vector fields $w_1,\dots,w_s$ commute to each other. It is also clear that 
$\tilde w_i$ are of type $(1,0)$ on $X_r$ and so are  the vector fields $w_i$.
It only remains to show that $w_i$ are in fact holomorphic on $X_r$.

Let $(t^1, \dots , t^s, z^1, \dots , z^{n-s})$ be local coordinates
fulfilling (\ref{local_coordinates}).
A local basis of vector fields of type $(1,0)$ on $X_r$ is given by
\begin{equation}\label{basis10}
\big\{\tilde w_1,\dots,\tilde w_s,\frac{\partial}{\partial z^1}-
\sum_{j=1}^s c^j_1 \overline{v}_j,\dots, \frac{\partial}{\partial z^{n-s}}-\sum_{j=1}^s c_{n-s}^j \overline{v}_j  \big\},
\end{equation}
where $c_l^j= \sum_{k=1}^{n-s}b^j_k\beta^k(\partial/{\partial z^l})$. Notice that the functions $c_l^j$
are holomorphic on $X_r$ and basic with respect to the foliation $\mathcal F$, that is $c_l^j$ only depend on the coordinates $z^1, \dots , z^{n-s}$. The vector fields $w_i$ are holomorphic if and only if 
the Lie brackets $[w^i,\overline w]$ are of type $(0,1)$ on $X_r$ for each local vector field $w$ of type $(1,0)$ 
on $X_r$. But this is a straightforward computation using the local basis of $T X_ r^{1,0}$ given above.
\end{proof}

\begin{remarks}\label{rem:t-def}
(a) Notice that the almost complex structure on $X$ by the decomposition 
$T^*X^{\mathbb C} = N_r^{1,0}\oplus\overline{N_r^{1,0}}$
is real-analytic. 

(b) The matrix $C$ in $r=(C,\theta)$ corresponds to a choice of a basis of
fundamental vector fields of the action.  Hence, if  $r=(C,\theta)$ and  $r'=(C',\theta)$ with $C'\neq C$ 
then the complex manifolds $X_r$ and $X_{r'}$ are identical, as well as the holomorphic foliations
$\mc F_r$ and $\mc F_{r'}$, although the actions  $\varpi_{r}$ and  $\varpi_{r'}$ are different.

(c) One can consider what, in principle, could be a more general class of tangential deformations of
the action, namely those determined by the global 1-forms $\tilde\alpha^\ell$ on $X$ defined as
\begin{equation}\label{mes_alpha}
\tilde\alpha^\ell =  \alpha^\ell + \sum \overline e_j^\ell\overline\alpha^j  + \sum \overline f_k^\ell \overline\beta^k + 
\sum g_j^\ell \alpha^j + \sum h_k^\ell \beta^k,
\end{equation}
where the coefficients $e_j^\ell, f_k^\ell, g_j^\ell, h_k^\ell$ are complex numbers. In fact, if these coefficients are
small enough then there are smooth vector fields uniquely defined by the conditions 
$\tilde\alpha^\ell(\tilde v_j)=\delta_j^\ell$ and 
$\overline{\tilde\alpha}^\ell(\tilde v_j) = 0$. Moreover, if we denote by $\tilde Q$ the real subbundle of $T^*X^{\mbb C}$ generated by
the 1-forms $\tilde\alpha^\ell$, then $N^{1,0} = \tilde Q\oplus\nu^*\mc F^{1,0}$ defines an integrable almost complex structure and the vector fields $\tilde v_j$ are holomorphic on the new complex manifold $\tilde X$.

Notice however that this complex structure and the vector fields $\tilde v_j$ are the same as the ones 
determined by the 1-forms 
$$
\alpha^\ell + \sum \overline e_j^\ell\overline\alpha^j  + \sum \overline f_k^\ell \overline\beta^k + 
\sum g_j^\ell \alpha^j.
$$
Furthermore,  $\tilde X$ endowed with the $\mbb C^s$-action 
defined by the vector fields $\tilde v_j$  can be identified with a suitable $(X_r,\varpi_r)$ in the family 
constructed above. More precisely, if we denote by $D$ the invertible matrix $(\delta_j^\ell + g_j^\ell)$
then $\tilde X$ coincides with $X_r$, where $r=(C,(\theta^1,\dots,\theta^s))$ is determined by 
$C=(c_m^i) = D^{-1}$ and 
$\overline\theta^i = \sum\overline a^i_j\overline\alpha^j + \sum\overline b_k^i\overline\beta^k$ with $\overline a^i_j 
= \sum c_m^i \overline e_j^m$
and $\overline b^i_j = \sum c_m^i \overline f_j^m$.
Therefore this construction does not provide more general deformations of $(X,\varpi)$. In fact it is proved in the 
following theorem that $(X_R,\varpi_R)$ is the versal family of tangential deformations of $(X,\varpi)$. In particular 
the family is complete and it contains all the small tangential deformations of $(X,\varpi)$ up to isomorphism.

(d) It follows from the above remark that the family $(X_R,\varpi_R)$ does not depend on the choice
of the 1-forms $\alpha^i$ fulfilling $\alpha^i(v_j)=\delta^i_j$. Therefore that family is naturally 
associated to $(X,\varpi)$.
\end{remarks}

\begin{remark} \label{rem:t-def_bis}
Let $\varpi'$ be a locally free $\mbb C^s$-action on a compact complex manifold $X'$ and assume
that $(X',\varpi')$ coincides with a pair $(X_r,\varpi_r)$ in the family $(X_R,\varpi_R)$. Then $(X,\varpi)$
coincides with an element of the family  $(X'_R,\varpi'_R)$ associated to $(X',\varpi')$.
Therefore the original action $\varpi:G\times X\rightarrow X$ is also a tangential deformation of
$\varpi_{r}:G\times X_{r}\rightarrow X_{r}$ for each $r\in R$.
\end{remark}

\begin{theorem}\label{teo:versal_space}
The family $(X_R,\varpi_R)$ defined above is the versal family of tangential deformations
of $(X,\varpi)$.
\end{theorem}

\begin{proof}
Let $\kappa\!: T_0R\to h^1(X,TX)$ denote the Kodaira-Spencer map associated to the family
$(X_R,\varpi_R)$ at $0$. Let $(C = {\mathrm id} + C',\theta)$ be a given element in $\Xi$ close
to $({\mathrm id}, 0)$. Then $(C',\theta)$ is an element of $T_0R = 
\operatorname{End}\mathbb C^s \oplus H^0(X,\Omega^1_X)^s$, close to zero, that can be seen as a cocycle in 
$C_h^1(\mbb C^s,\Omega^{0,0}(X,T\mathcal F))\oplus C_h^0(\mbb C^s,\Omega^{0,1}(X,T\mathcal F))$.
An easy computation shows that $\kappa(C',\theta)$ is just the cohomology class in $H^1_{\mbb C^s}(X, T\mc F)$
of that cocycle. This implies that: (i) the restriction of $\kappa$ to 
$\operatorname{End}\mbb C^s$ is just the map $\tau$ in~(\ref{eq:AbelianShortSequence}),
and (ii) the restriction of $\kappa$ to $T_0R_0\equiv H^0(X,\Omega^1_X)^s$ is a section of the morphism $\chi$
in~(\ref{eq:AbelianShortSequence}). Hence $\kappa$ is an isomorphism and the statement follows 
from Proposition~\ref{prop:criteri_versal}.
\end{proof}

Let $\K_X$ denote the Kuranishi space of $X$, that is the parameter space of the versal family of 
deformations of the complex manifold $X$. Recall that $\K_X$ is a (possibly singular) analytic space 
whose tangent space at the distinguished point $0\in \K_X$ is naturally identified to $H^1(X,\Theta_X)$.
The restriction $X_{R_0}$ to $R_0$ of the family $X_R$ can be seen as a family of deformations 
of the complex structure of $X$. Hence there is a well defined forgetful map of germs of analytic spaces
$$
\phi\!: (R_0, 0) \longrightarrow (\K_X,0).
$$ 
Moreover its tangent map $d_0\phi$ at $0$ coincides with the morphism 
$$H^1(X,T\mc F)\to H^1(X,\Theta_X)$$
induced by the decomposition $TX = T\mc F \oplus T\mc G$. More precisely, given an element 
$\theta = (\theta^{1}, \dots,\theta^{s})$ in $T_0R_0$ with
$\theta^{i} = \sum a_{j}^{i} \alpha^{j} + \sum b_{k}^{i}\beta^{k}$, its image $d_0\phi(\theta)$
is the vector valued $(0,1)$-form
$$
\tilde\theta = \sum \bar\theta^{i}\otimes v_i.
$$
Since holomorphic forms are closed and non exact, one can see using local coordinates  
fufilling (\ref{local_coordinates}) that $\tilde\theta$
defines a cohomology class in $H^1(X,\Theta_X)$ which is not trivial unless $\theta^i=0$ for all $i=1,\dots s$.
This implies that $\phi$ is an embedding proving the following statement.

\begin{theorem}\label{teo:kuranishi}
Let $X$ be a compact manifold in the class $\mc C$ such that 
$s=\dim_{\mbb C} T_X= \dim_{\mbb C}\mf h_X/\mf h^1_X\neq 0$. 
Then the Kuranishi space $\K_X$ of $X$ contains a 
smooth subspace of dimension $s\cdot b_1(X)/2$.
More precisely, the canonical morphism 
$H^1(X,T\mc F)\equiv H^1(X,\mc O_X)^s\to H^1(X,\Theta_X)$ is injective and its image are 
unobstructed infinitesimal deformations of the complex structure of $X$.
\end{theorem}
\medskip

We shall make distinction of two particular types of tangential deformations defined as follows
\begin{definition}
Let $r = (C, \theta)\in R$ be given and set
$\theta =(\theta^{1},\dots,\theta^{s})$.
We say that the associated tangential deformation $(X_{r}, \varpi_{r})$ is
\begin{enumerate}[(i)]
\item a {\sl purely tangential deformation}, or $t$-deformation, of $(X, \varpi)$ if
$\theta^{i}\in \langle\alpha^{1},\dots, \alpha^{s}\rangle$, i.e. if the coefficients
$b_{k}^{i}$ in (\ref{new-alpha}) are all zero,
\item a {\sl basic deformation}, or $b$-deformation, of $(X, \varpi)$ if
$\theta^{i}\in \langle\beta^{1},\dots, \beta^{p-s}\rangle$, i.e. if the coefficients $a_{j}^{i}$
in (\ref{new-alpha}) are all zero.
\end{enumerate}
\end{definition}

\begin{remarks}\label{rem:ComplFol}
(a) The holomorphic and invariant foliation $\mc G$, which is transversal and complementary
to $\mc F$, remains holomorphic and invariant for each $t$-deformation 
of $(X,\varpi)$. On the other hand the foliation $\mc F$ remains holomorphic for each $b$-deformation
of $(X,\varpi)$.

(b) If $(X_r,\varpi_r)$ is a $t$-deformation (resp. $b$-deformation) of $(X,\varpi)$ then $(X,\varpi)$
is a $t$-deformation (resp. $b$-deformation) of $(X_r,\varpi_r)$. 

Remark also that  each small tangential 
deformation $(X_{r'},\varpi_{r'})$ of $(X,\varpi)$ can be seen as a $b$-deformation 
of a $t$-deformation $(X_r,\varpi_r)$ of $(X,\varpi)$. More precisely, if $X_{r'}$ is determined 
by 1-forms $\tilde{\alpha}^i=\alpha^i+\sum \overline{a}^i_j \overline{\alpha}^j+
\sum \overline{b}^i_k \overline{\beta}^k$ then the manifold $X_r$ determined by the 1-forms 
${\hat\alpha}^i=\alpha^i+\sum \overline{a}^i_j \overline{\alpha}^j$ is a $t$-deformation of $X$. Moreover, since 
$\tilde{\alpha}^i={\hat\alpha}^i+\sum \overline{b}^i_k \overline{\beta}^k$ the manifold $X_{r'}$ is a $b$-deformation 
of $X_r$. 

Notice that, as $X$ is a tangential deformation of $X_r$, the previous argument also proves 
that $(X_{r'},\varpi_{r'})$ can be seen as a $t$-deformation of a $b$-deformation of $(X,\varpi)$. 
\end {remarks}
\smallskip

As it is remarked in the proof of Proposition~\ref{prop:integrabilitat}, local holomorphic functions on $X$
that are basic with respect to $\mathcal F$ are also holomorphic on $X_r$ and basic with respect to $\mc F_r$.
By definition, the global 1-forms $\alpha_{r}^{i}$ are of type $(1,0)$
on $X_{r}$ and they are also closed by construction. Hence $\alpha_{r}^{i}$
are holomorphic forms on $X_{r}$. It follows that, for each $r\in R$,  one has
$h^{1,0}(X_r) = h^{1,0}(X)$,  $h^{0,1}(X_r) = h^{0,1}(X)$ and
$b_1(X) = h^{1,0}(X_r) + h^{0,1}(X_r)$. In particular, the space of holomorphic 1-forms
on $X_r$ can be decomposed as
\begin{equation}\label{eq:espai-formes}
H^0(X_r,\Omega^1_{X_r})=
\langle \alpha_r^1,\dots,\alpha_r^s\rangle \oplus \langle \beta^1,\dots,\beta^{p-s} \rangle.
\end{equation}

\begin{proposition}\label{prop:mateixos_camps}
Let $\mf a_X$ be an Abelian algebra fulfilling ${\mf h}_X = \mf a_X\oplus {\mf h}^1_{X}$
and let $X_r$ be a tangential deformation of $X$. 
\begin{enumerate}[\rm (i)]
\item If $\mf a_X$ is in the center of ${\mf h}_X$ then ${\mf h}^1_{X_r}={\mf h}^1_{X}$.
\item If $X_r$ is a $b$-deformation of $X$ then $\mf a_X \subset \mf h_{X_r} $ and 
$\mf h_{X_r} =  \mf a_X\oplus \mf h^1_{X_r}$.
\end{enumerate}
If both conditions are fulfilled then $\mf h_{X_r} = \mf h_X$.
\end{proposition}

\begin{proof}
(i) Let $w\in\mf h^1_X$ be given. 
Notice that  $\alpha^i(w)=\alpha^i_r(w)=0$ for $i=1,\dots,s$. Since $[\mf a_X, \mf h^1_X]=0$,
$w$ is written locally as
$$w=\sum_{k=1}^{n-s} f^k (z) \frac{\partial}{\partial z^k},$$
where $f^k$ are holomorphic functions on $z^{1},\dots,z^{n-s}$. It follows that
$w$ is of type $(1,0)$ on $X_r$. Now the same argument used in the proof of 
Proposition~\ref{prop:integrabilitat} shows that $w$ is in fact holomorphic in $X_r$.
The original manifold $X$ is a tangential deformation of $X_r$ (cf. Remark~\ref{rem:t-def_bis}), 
so $\mf h^1_{X_r}$ 
does not contain any additional tangent fields besides those arising in $X$.

(ii) If $X_r$ is a $b$-deformation of $X$ then the elements $v\in \mf a_X$ are still 
holomorphic on $X_r$, as it follows from the identities (\ref{new-vectors}) to 
(\ref{new-vectors-last}).
\end{proof}
\medskip

We discuss now the tangential deformations of the examples of K\"{a}hler manifolds 
introduced in Section~\ref{Sec:examples}.

\begin{example}[Deformations of tori]
If $X$ is a complex torus $\mbb T^s = \mbb C^s/\Lambda$ then $H^0( \mbb T^s,\Omega^1_{\mbb T^s})^s$
can be identified with the dual space of $H^1(\mbb T^s, \Theta_{\mbb T^s})$. 
In that case the space $R$ is just the product of
$\operatorname{GL}(s,\mbb C)$ with the versal space of deformations of the complex manifold
$\mbb T^s$ as it is described by Kodaira and Spencer in \cite{KodSpe}.
\end{example}

\begin{example}[Deformations of principal torus bundles]
Assume that $X$ is the total space of a principal torus bundle
$\mbb T^s \rightarrow X \rightarrow B$ where $\mbb T^s=\mbb C^s/\Lambda$.
As it was already noticed  (cf. Example~\ref{example_fibration}), when $X$ is K\"ahler
then $B$ is also a
K\"ahler manifold and the bundle is flat. The possible choices of the $s$-tuple
$(\alpha^1, \dots, \alpha^s)$ correspond exactly to the different flat connections on that bundle.

In this situation the space of basic deformations of the action is naturally identified to
$H^0(B, \Omega^1_B)^s \cong H^1(B, \mc O^s_B)$. The exact sequence of sheaves
$
0\rightarrow \Lambda \rightarrow \mc O^s_B \rightarrow \mc O^s_B/\Lambda \rightarrow 0
$
over $B$ induces the exact cohomology sequence
$$
\dots \longrightarrow H^1(B, \Lambda) \longrightarrow H^1(B, \mc O_B^s)
\longrightarrow  H^1(B, \mc O_B^s/\Lambda) \stackrel{\tau}{\longrightarrow} H^2(B, \Lambda)
\longrightarrow \dots
$$
It follows that $H^1(B, \mc O_B^s)$ is maped onto $\ker \tau$, which is the space of
(isomorphism classes of) topologically trivial $\mbb C^s/\Lambda$-principal bundles over $B$. 
Replacing $B$ by a finite (unramified) covering we can suppose that $H^2(B, \Lambda)$
has no torsion. In that case each
topologically trivial principal torus bundle over $B$ can be obtained from $X$ by a suitable
basic deformation of the action and, in particular, the trivial bundle is a tangential deformation of
$(X,\varpi)$.

In general, each tangential deformation of $(X,\varpi)$ is a principal torus bundle over $B$.
Purely tangential deformations  $r\in R$ just deform the fiber $\mbb C^s/\Lambda_r$.
\end{example}

\begin{example}[Deformations of suspensions]\label{ex:def_suspensions}
Assume that $X$ is a suspension over a torus $\mbb T^s = \mbb C^s/\Lambda$, 
that is $X = F\times_\Lambda \mbb C^s$. It follows from Remark~\ref{rem:ComplFol}~(a)
that each purely tangential deformation of $X$ is still a
suspension, although this is not true for general tangential deformations of $X$.

If $X$ is a suspension one has $H^0(X,\Omega^1_{X/\mc F}) = H^0(F,\Omega^1_F$). Therefore if 
$b_1(F)=0$, or what is equivalent $b_1(X)=2s$, then $X$ has no basic deformations
and each tangential deformation of $X$ is again a suspension.
\end{example}

We end this section stating several general properties of tangential deformations of compact manifolds 
of the class $\C$ endowed with a locally free $\mbb C^s$-action.
\medskip

If $\tilde X \rightarrow X$ is an unramified covering then $\tilde X$ is also endowed with 
a locally free $\mbb C^s$-action in a natural way. Notice however that the maximal rank of 
such an action on $\tilde X$ can increase, that is $\dim_{\mbb C} \mf h_{\tilde X}/\mf h^1_{\tilde X} \geq s$.
The following statement follows straightforward from
the definitions.

\begin{proposition}\label{prop:finite_covering}
A finite covering $X'_r$ of a tangential deformation $X_r$ of $X$ is a tangential deformation
of a finite covering $X'$ of $X$. Conversely, a tangential deformation $X'_r$ of a finite covering 
$X'$ of $X$ is a finite covering of a tangential deformation $X_r$ of $X$.
\end{proposition}

\begin{proposition}\label{prop:ruled_estable}
Let $X$ be a geometrically ruled manifold in the class $\mc C$. Then $X_r$ is a geometrically 
ruled compact complex manifold for each $r\in R$.
\end{proposition}

\begin{proof}
Let $X$ be a compact manifold in the class $\mc C$ endowed with a locally free 
$\mbb C^s$-action $\varpi$, and assume that $X$ is the total space of a locally trivial 
fiber bundle, $\varphi\!:X \to B$, with fiber 
$\mbb C P^k$ and defined by a cocyle $\{g_{ij}\}$ in $H^1(B,\operatorname{PGL}(k+1, \mc O_B))$.
Then $B$ is also in the class $\mc C$. Since automorphisms of $X$ close to the identity 
must preserve the fibration (cf. \cite{Bla}) each fundamental vector field $v$ of the action 
$\varpi$ projects into a non vanishing vector field $\hat v$ on $B$. These vector fields 
define a locally free $\mbb C^s$-action $\hat\varpi$ on $B$ and the projection $\varphi$ is 
equivariant with respect to the $\mbb C^s$-actions. 
Notice that every global holomorphic 1-form on $X$ vanishes on the fibers of $\varphi$. Hence 
$\varphi^*\!:H^0(B,\Omega^1_B)\to H^0(X,\Omega^1_X)$ is an isomorphism. Using this 
identification, one can associate, to each tangential deformation $X_r$ given by a pair $r =(C, \theta)$, 
a tangential deformation $B_r$ defined by the same expression~(\ref{new-alpha}).
Then the projection $\varphi\!:X_r\to B_r$ is  holomorphic. 

Now, as the $\mbb C^s$-action $\varpi$ is transverse and preserves the fibration, the cocycle 
$\{g_{ij}\}$ can be chosen to be constant along the orbits of the $\mbb C^s$-action $\hat\varpi$ on $B$. 
This implies that the matrix-valued functions $g_{ij}$ are still holomorphic on $B_r$ and therefore they define
$X_r$ as a geometrically ruled manifold over $B_r$.
\end{proof}

\begin{proposition}
A tangential deformation $X_r$ of $X$
uniquely determines a tangential deformation $\A{X}_r$
of the Albanese torus $\A{X}$ of X such that $\A{X}_r = \A{X_r}$,
i.e. $\A{X}_r$ is the Albanese torus of  $X_r$.
\end{proposition}

\begin{proof}
The locally free $\mbb C^s$-action $\varpi$ on $X$ induces a locally free $\mbb C^s$-action 
$\hat\varpi$ on the corresponding Albanese torus $\A{X}$ and the natural map $\phi\!: X\to\A{X}$ is 
equivariant with respect to these actions. Since the map $\phi$ induces an 
isomorphism $\phi^*\!:H^0(\A{X}, \Omega^1_{\A{X}})\to H^0(X, \Omega^1_X)$, there is 
a natural correspondence, as the one considered in the proof of the above proposition, 
associating to each tangential deformation $X_r$ of $X$ a tangential deformation 
$\A{X}_r$ of $\A{X}$. It follows from the definitions that $\A{X}_r=\A{X_r}$.
\end{proof}

\begin{proposition}\label{prop:stability_dimkod}
The manifolds $X$ and $X_r$ have the same plurigenera for all $r\in R$. 
In particular the Kodaira dimension $\operatorname{kod}(X_{r})$ of $X_{r}$ is equal to
$\operatorname{kod}(X)$ for each $r\in R$.
\end{proposition}

\begin{proof}
Let $(t^1,\dots,t^s, z^1,\dots,z^{n-s})$ be local coordinates of $X$ 
fulfilling~(\ref{local_coordinates}) and let $\alpha_r^i$ be the global 
holomorphic 1-forms on $X$ defined by~(\ref{new-alpha}). Notice that $z^j$ 
are also holomorphic functions on $X_r$. A local section of $K_{X_r}^{\otimes m}$
is of the form $f(t,z) \eta^m$ where 
$\eta=\alpha^1\wedge\dots\wedge\alpha^s\wedge dz_{i}^1\wedge\dots\wedge dz_{i}^{n-s}$
and $f=f(t,z)$ is a holomorphic function.

Nakamura and Ueno have proved that, for each compact complex manifold $Y$, the group 
$\operatorname{Aut}^0_{\mbb C}(Y)$ acts trivially on $H^0(Y,K_Y^m)$ (cf. \cite[Corollary 14.8]{Ue}). 
Applying this result to the $\mbb C^s$-action we see that $f$ is 
necessarily a basic function, i.e. $f=f(z)$ only depends on the coordinates
$(z^1,\dots,z^{n-s})$ transverse to $\mc F$. Therefore $H^0(X_r,K_{X_r}^{\otimes m})$
is naturally identified to the space of basic global sections of the line bundle
$(\det \nu^*\mc F_r)^m$
and this space is independent of $r\in R$.
\end{proof} 

If $X$ is a K\"{a}hler manifold then Kodaira's stability theorem implies that $X_r$ is also 
K\"{a}hler for every $r$ in a neighborhood of $(\mathrm{id},0)$ in $R$, and we will prove the following 
result, stating that the class of K\"{a}hler manifolds is stable under tangential deformations. 
In sharp contrast the class $\mc C$ is not stable under 
small deformations. We do not know if tangential deformations of Fujiki manifolds 
still belong to the class $\mc C$.

\begin{theorem}\label{teo:kahler_stability}
Assume that $X$ is a K\"{a}hler manifold, then $X_r$ are also K\"{a}hler manifolds for 
each $r\in R$.
\end{theorem}

We postpone the proof of this theorem to Section~\ref{sec:structure}, where the assertion will 
be a corollary of general results on the structure of such manifolds. 
Nevertheless, accepting that the theorem has already 
been proved we can state the following two propositions:

\begin{proposition} \label{prop:hodge_numbers}
Assume that $X$ and $X_r$ belong to the class $\mc C$. Then 
\begin{enumerate}[\rm (i)]
\item $h^{k,0}(X) = h^{k,0}(X_{r})$ and $h^{0,k}(X) = h^{0,k}(X_{r})$ for each $k$,
\item $h^{1,1}(X) = h^{1,1}(X_{r})$. \end{enumerate}
If $X$ is a K\"{a}hler manifold the above identities hold for each $r\in R$.
\end{proposition}

\begin{proof} 
Lemma~\ref{decomposition} implies that there is a well defined $\mbb C$-linear
map from $H^0(X, \Omega^{k}_X)$ to $H^0(X_{r}, \Omega^{k}_{X_r})$ determined 
by the condition
$$
\gamma = \sum_{0\leq |I|\leq \ell} \alpha^{I} \wedge \gamma^I  \longmapsto
\gamma_{r} = \sum_{0\leq |I|\leq \ell} \alpha^{I}_{r} \wedge \gamma^I.
$$
where $\ell =  \min(k,s)$,
$ \alpha^{I}_{r} = \alpha^{i_1}_{r} \wedge \dots \wedge \alpha^{i_{|I|}}_{r}$, 
$\alpha_{r}^{i}= \alpha^{i} + \overline\theta^{i}$ and $\gamma^I$ are global basic forms.
Clearly this map is in fact an isomorphism. Preservation of $b_2, h^{2,0}, h^{0,2}$
implies that of $h^{1,1}$.
\end{proof}

\begin{proposition}\label{prop:tot_versal}
Assume that $X$ is a K\"{a}hler manifold. Then the family $(X_R,\varpi_R)$
is the versal family of $(X_r,\varpi_r)$ at each point $r\in R$.
\end{proposition}

\begin{proof}
In that case $X_r$ is a K\"{a}hler manifold for each $r\in R$ and 
$H^1(X_r, \mathcal O_{X_r})$ is isomorphic to $H^0(X_r,\Omega^1_{X_r})$. Moreover,
the dimension of these vector spaces does not change with $r\in R$. Hence the argument 
used to prove Theorem~\ref{teo:versal_space} also applies to each $X_r$ with $r\in R$.
\end{proof}

%%%%%%%%%%%%%%%%%%% SECTION 6

\section{The approximation theorem}\label{sec:approximation}

In this section we prove a key result for describing the structure of  compact complex manifolds $X$
in the class $\mc C$ endowed with a locally free holomorphic $\mbb C^s$-action. Namely, we prove that
such a manifold can be approximated 
by tangential deformations $X_\epsilon$ of $X$ that are suspensions over a 
torus of dimension $s$, i.e. $X_{\epsilon}$ belonging to the class of manifolds described in 
Example~\ref{example_suspension}. 

We shall make use of the following criterium, which is valid for
arbitrary compact complex manifolds.

\begin{lemma}\label{CriteriSusp} Let $X$ be a compact complex manifold endowed with a
locally free action of $\mbb C^s$ and let $v_1,\dots,v_s$ be a basis of fundamental vector fields 
of the action. Then $X$ is a suspension over a torus $T= \mbb C^s/\Lambda$, i.e. 
$X = F\times_\Lambda \mbb C^s$ for a suitable representation 
$\Lambda\to \operatorname{Aut}_{\mbb C}(F)$,
if and only if one can choose holomorphic 1-forms 
$\alpha^j$ on $X$ fulfilling $\alpha^j(v_i)=\delta^j_i$ and such that
the set of periods of $\alpha^1,\dots,\alpha^s$ 
$$
\Lambda = \{ (\int_\gamma \alpha^1,\dots , \int_\gamma \alpha^s) \mid \gamma\in
H_1(X,\mbb Z) \}
$$
is a lattice of $\mbb C^s$. \end{lemma}

\begin{proof}
For each $i = 1, \dots, s$ let $f^i$ be a multivalued holomorphic function such that
$\alpha^i = df^i$. If the set of periods $\Lambda$ is a lattice of $\mbb C^s$ then
$f=(f^1, \dots, f^s)$ is a well-defined holomorphic submersion from $X$ onto
$\mbb C^s/\Lambda$. The converse is clear.
\end{proof}

\begin{remark} \label{rem:fibres_connexes}
The fiber space $F$ in the above theorem can always be taken connected. 
Indeed, if the fibers of the submersion $f : X \rightarrow T$ are not connected,
one can consider the Stein factorization of  $f$
$$
\begin{diagram}
\node{X} \arrow{e,t}{p} \arrow{se,r}{f} \node{N} \arrow{s,r}{\pi} \\
\node{}\node{T}
\end{diagram}
$$
Since the map $p$ is finite and $f$ is a submersion, the map $\pi$ is necessarily 
an unbranched cover, thus $N$ is also a torus that 
covers $T$. The map $p$ is then the sought suspension.
\end{remark}

The following result is an improvement of Theorem 2.16 in \cite{Man}. For our purposes here, 
an essential point of the theorem is the fact that the manifold $X_\epsilon$, which approximates $X$ and 
is a suspension over a torus, can be taken to be a basic deformation of $X$. 

\begin{theorem}\label{thm:Rio}
Let $X$ be a compact complex manifold in the class $\mc C$ endowed with a
locally free $\mbb C^s$-action. There is an arbitrarily small $b$-deformation
$X_{\epsilon}$ of $X$ which is a suspension over a complex torus $T= \mbb C^s/\Lambda$.
More precisely, there is a connected compact complex manifold $F$ and a group representation
$\rho:\pi_1(T)=\Lambda\rightarrow \operatorname{Aut}_{\mbb C} (F)$ such that $X_{\epsilon}$
is the suspension of $\rho$.
\end{theorem}

\begin{proof}
First, we claim that there exist arbitrarily small holomorphic 1-forms 
$\theta^i,\eta^i$ on $X$ such that, if we set $\tilde\alpha^i= \alpha^i+\overline{\theta}^i+\eta^i$, then
$$
\Lambda=\big\{ \big(\int_\gamma \tilde\alpha^1,\dots , \int_\gamma \tilde\alpha^s\big)\mid\gamma\in
H_1(X,\mbb Z) \big\}
$$
is a lattice of $\mbb C^s$. The proof is essentially the same as the one given in \cite[p. 490]{Man}. 
Let $\{ \alpha^{1},...,\alpha^{s},\beta^{1},...,\beta^{p-s} \}$ be the basis of $H^0(X,\Omega^1_X)$
considered in~(\ref{eq:base_formes}). We decompose
$$
\theta^i= \theta^i_t + \theta^i_b =\sum a^i_j \alpha^j+\sum b^i_k \beta^k, 
\qquad \eta^i= \eta^i_t + \eta^i_b = \sum c^i_j \alpha^j + \sum d^i_k \beta^k.
$$

It follows from Remark~\ref{rem:t-def} (c) that the 1-forms $\alpha^i + \overline\theta^i + \eta_t^i$ determine a
complex manifold endowed with a $\mbb C^s$-action which is isomorphic to a suitable 
$(X_{r'},\varpi_{r'})$ close to 
$(X,\varpi)$. Since $\tilde\alpha^i$ are holomorphic 1-forms on $X_{r'}$, Lemma~\ref{CriteriSusp} 
implies that the manifold $X_{r'}$ is a suspension over a torus.   

As it has been noticed in Remark~\ref{rem:ComplFol}~(b), we can see  $(X_{r'},\varpi_{r'})$
as a $t$-deformation of a suitable $(X_r,\varpi_r)$ close to $(X,\varpi)$ with the property that 
$(X_r,\varpi_r)$ is a  $b$-deformation of $(X,\varpi)$. But $X_r$ is a suspension over a torus as 
$t$-deformations of suspensions are again suspensions (cf. Example~\ref{ex:def_suspensions}). 
Hence $(X_r,\varpi_r)$ fulfills the required conditions as Remark~\ref{rem:fibres_connexes}
assures that the fibre $F$ of the suspension can always be 
taken connected. 
\end{proof}

It follows from the above theorem that there is an exact sequence of groups
$
1\to\pi_1(F)\to \pi_1(X) \to \mbb Z^{2s}\to 1.
$

\begin{remarks}\label{rem:canvi_base}
(a) It follows from Remark~\ref{rem:t-def}~(a) that $X_\epsilon$ is real-analitically 
isomorphic to the manifold $X$. More precisely, as the two complex manifolds 
$X$ and $X_\epsilon$ have the same underlying real-analytic structure, the identity
map ${\mathrm id}\!: X \to X_\epsilon$ is real-analytic.

(b) Assume that $X$ is a suspension over a torus $T=\mbb C^s/\Lambda $, i.e. $X= F\times_\Lambda\mbb C^s$.
If $T'$ is a torus isogenous to $T$ then $T'= \mbb C^s/\Lambda'$, where $\Lambda'$ is a sublattice of
$\Lambda$, and there is an Abelian finite covering $X'$ of $X$ which is a suspension over $T'$. 
From this fact one deduces easily that if $X$ is a suspension over an $s$-dimensional torus, where
$s=\dim_{\mbb C}T_X = \dim_{\mbb C} \mf h_X/\mf h^1_X$, then $X$ has 
an Abelian finite covering  $X'$ that is a suspension over $T_X = \Phi( \operatorname{Aut}^0_{\mbb C} (X))
\subset \A{X}$.

(c) Let $X$ be a K\"{a}hler manifold in the hypothesis of the above theorem. Then 
$X_{\epsilon}$ is the suspension $F\times_\Lambda\mbb C^s$ of a representation
$\rho:\Lambda\to \operatorname{Aut}_{\mbb C} (F)$. 
Since $X_{\epsilon}$ is a K\"{a}hler manifold, there is a sublattice $\Lambda'$ of 
$\Lambda$ such that  $\rho(\Lambda')\subset \operatorname{Aut}^0_{\mbb C} (F)$
(cf. \ref{example_suspension}). Then the suspension 
$X'_{\epsilon}= F\times_{\Lambda'}\mbb C^s$ is an Abelian finite covering of  
$X_{\epsilon}$. Notice that $X'_{\epsilon}$ is topologically isomorphic to the 
product $F\times T'$, where $T' = \mbb C^s/\Lambda'$, and that $X'_{\epsilon}$
can be seen as a $b$-deformation of a finite covering $X'$ of $X$
(cf. \ref{prop:finite_covering} and \ref{rem:t-def_bis}).

(d) Any (unramified) finite covering $X'$ of $X$ is naturally endowed with a $\mbb C^s$-action
but the maximal rank can increase, that is 
$$\dim_{\mbb C} \mf h_{X'}/\mf h^1_{X'} \geq \dim_{\mbb C} \mf h_{X}/\mf h^1_{X}.$$
Examples where this inequality is strict can be constructed by considering suspensions of 
bielliptic surfaces.
\end{remarks}

A more precise description of the fiber manifold $F$ and of the representation 
$\rho:\Lambda\to \operatorname{Aut}_{\mbb C} (F)$ in the above theorem 
can be obtained by considering (unramified) finite coverings $X'$ of the manifold $X$.
This is the content of the following proposition. It can be seen as an improvement of the 
approximation theorem and will play a key role in the structure theorems of the next Section.

\begin{proposition}\label{prop:primera_reduccio}
Let $X$ be a compact K\"{a}hler manifold such that 
$s=\dim_{\mbb C}T_X = \dim_{\mbb C}\mf h_X/\mf h^1_X >0$.
There is a finite Abelian covering $X'_\epsilon$
of an arbitrarily small $b$-deformation $X_\epsilon$ of $X$, which is the suspension
over a torus $T= \mbb C^{s'}/\Lambda$, of dimension $s'\geq s$, 
associated to a representation 
$\rho\!: \Lambda\to \operatorname{Aut}_{\mbb C}(F)$
with the following properties:
\begin{enumerate}[\rm (i)]
\item $F$ is a connected compact K\"{a}hler manifold with $\mf h_{F}/{\mf h^{1}_{F}} = 0$, i.e.
$F$ has no non-singular holomorphic vector fields,
\item $\rho(\Lambda) \subset \operatorname{Aut}^0_{\mbb C}(F)$.
\end{enumerate}
\end{proposition}

\begin{remark}
The manifold $X'_\epsilon$ can also be seen as a (small) $b$-deformation of a finite
Abelian covering $X'$ of $X$. This follows from Proposition~\ref{prop:finite_covering}. 
\end{remark}

\begin{proof}
Let $\mf a_X$ be an central Abelian subalgebra of $\mf h_X$ 
fulfilling $\mf h_X= \mf a_X \oplus \mf h^1_X$. 
It defines a locally free 
$\mbb C^s$-action on $X$ and, by the above Theorem~\ref{thm:Rio} there is an 
arbitrarily small $b$-deformation $X_\epsilon$ of $X$ which is the suspension
$F_1\times_{\Lambda_1} \mbb C^s$ associated to a representation 
$\rho_1\!:\Lambda_1\to\operatorname{Aut}_{\mbb C}(F_1)$.
Notice here that, applying Proposition~\ref{prop:mateixos_camps}, 
we have $\mf a_{X} \subset \mf h_{X_\epsilon}$. Therefore we can chose 
$\mf a_{X_\epsilon} = \mf a_X$ as the complementary of $\mf h^1_{X_\epsilon}$
in $\mf h_{X_\epsilon}$.

As in Remark \ref{rem:canvi_base} (c) above, we consider 
the sublattice $\Lambda'_1 = \rho_1^{-1}(\operatorname{Aut}^0_{\mbb C}(F_1))$ of 
$\Lambda_1$. 
Then the suspension 
$X^1_{\epsilon}= F_1\times_{\Lambda'_1}\mbb C^s$ is an Abelian finite covering of  
$X_{\epsilon} = X^1_{\epsilon}/\Gamma_1$ where $\Gamma_1 = \Lambda_1/\Lambda'_1$. 
More precisely, each vector field on $X_{\epsilon}$ lifts to a unique vector field
on $X^1_{\epsilon}$ and therefore we can identify $\mf a_{X_\epsilon} = \mf a_X$ to an Abelian 
subalgebra of $\mf h_{X^1_\epsilon}$ that we still denote by $ \mf a_X$. Then there are 
$v_i\in  \mf a_X$ such that the monodromies $\rho_1(\gamma)$ for the set of elements $\gamma \in \Lambda_1$
generating $\Gamma_1$ lie in $\langle \gamma_1,\dots,\gamma_k \rangle$, where
$\gamma_i = \exp_1(v_i)$. 

If $F_1$ has no vector fields without zeros we are done. So, assume that 
$\mf h_{F_1}/{\mf h^{1}_{F_1}} \neq 0$. In that case we choose a subalgebra $\mf a_{F_1}$ of 
the center of $\mf h_{F_1}$ such that $\mf h_{F_1} = \mf a_{F_1} \oplus \mf h^1_{F_1}$
(cf. Proposition \ref{prop:touzet}). 
Since $\rho_1(\Lambda'_1)$ is contained in $\operatorname{Aut}^0_{\mbb C}(F_1)$ 
each element $\hat w$ of  $\mf a_{F_1}$  is $\rho_1(\Lambda'_1)$-invariant. 
This implies that 
$\hat w$ is the restriction of a globally defined holomorphic vector field $w$ on 
$X^1_\epsilon$ without zeros and commuting to $\mf a_{X}$. More precisely, 
$\hat w$ is the projection over  $F_1\times_{\Lambda_1}\mbb C^s$ of the vector field $(w,0)$
on the product manifold $F_1\times\mbb C^s$. Remark that, by construction, 
$\mf a_{X^1_\epsilon} = \mf a_{X}\oplus \mf a_{F_1}$ is Abelian and that 
$\mf h_{X^1_\epsilon} = \mf a_{X^1_\epsilon} \oplus \mf h^1_{X^1_\epsilon}$.
Denote $r=\dim \mf a_{X^1_\epsilon} > s$, and repeat the above construction 
beginning with $X^1_\epsilon$ instead of $X$. 

There is a small $b$-deformation $X^1_\delta$ of $X^1_\epsilon$ which is the
suspension $F_2\times_{\Lambda_2}\mbb C^r$ over a torus $\mbb C^r/\Lambda_2$,
associated to a representation 
$\rho_2\!:\Lambda_2\to\operatorname{Aut}_{\mbb C}(F_2)$. 
We see, using Proposition~\ref{prop:finite_covering}, that 
$X^1_\delta$ is a finite Abelian covering of a $b$-deformation $X_\delta$ of $X$. 
And, as above, $\mf a_{X^1_\epsilon}$ is an Abelian subalgebra of $\mf h_{X^1_\delta}$
with $\mf h_{X^1_\delta} = \mf a_{X^1_\epsilon} \oplus \mf h^1_{X^1_\delta}$.
In particular, $\gamma_i = \exp_1(v_i)$ can be seen as holomorphic transformations of
$X^1_\delta$.

We set $\Lambda'_2 = \rho_2^{-1}(\operatorname{Aut}^0_{\mbb C}(F_2))$ and
we define $X^2_\delta = F_2\times_{\Lambda'_2}\mbb C^r$. Then $X^2_\delta$
is a finite Abelian covering of $X^1_\delta = X^2_\delta/\Gamma_2$ where
$\Gamma_2 = \Lambda_2/\Lambda'_2$. 
Also as above, $\mf a_{X^1_\epsilon} = \mf a_{X}\oplus \mf a_{F_1}$ is naturally
included in $\mf h_{X^2_\delta}$. It follows in particular that the automorphisms 
$\gamma_i = \exp_1(v_i)$ for $v_i \in \mf a_X$ can be lifted as biholomorphisms $\tilde\gamma_i$ 
of $X^2_\delta$. 
This already implies that the composition 
$$X^2_\delta \longrightarrow X^1_\delta \longrightarrow X_\delta
$$
is a regular covering. 
Moreover, we can realize the monodromies $\rho_2(\gamma)$ for a set of elements in $\Lambda_2$
generating $\Gamma_2$ within a subgroup of automorphisms $\langle \hat\gamma_1,\dots,\hat\gamma_m \rangle$, where 
$\hat\gamma_j=\exp_1(w_j)$ and $w_j\in \mf a_{X^1_\epsilon}$. Since all the 
biholomorphisms $\tilde\gamma_i$ and $\hat\gamma_j$ are exponentials of vector fields
belonging to the same Abelian algebra, they commute to each other proving that 
the regular covering $X^2_\delta \to X_\delta$ is in fact Abelian. 

If $F_2$ does not have non-singular vector fields we have done. If this is not the case we 
repeat again the above construction starting with $X^2_\delta$.
Since the rank of 
locally free actions is bounded by the dimension $n$ of $X$ this procedure ends 
after a finite number of steps. 
\end{proof}

\begin{remark}\label{rem:encaix}
Let $\mf h_X= \mf a_X \oplus \mf h^1_X$, with $\mf a_X$ central, and 
$X'_\epsilon = F\times_\Lambda \mbb C^{s'}$
be as in the above Proposition. Along the proof it is shown that the Lie algebra decomposition 
$\mf h_{X'_\epsilon}= \mf a_{X'_\epsilon} \oplus \mf h^1_{X'_\epsilon}$ can be taken with 
the property that $\mf a_{X'_\epsilon}$ is generated by the projection over $X'_\epsilon$ 
of linear vector fields on $\mbb C^{s'}$ and that the lift to $X'_\epsilon$ of $\mf a_X= \mf a_{X_\epsilon}$ is
included in $\mf a_{X'_\epsilon}$.
\end{remark}

We are ready now to establish our refinement of Calabi's Theorem~\ref{thm:calabi2}, which 
was originally stated by Calabi in~\cite{Cal} with solvable Galois group $\Gamma$,. It was 
proved by Bogomolov in \cite{Bogo} for complex projective manifold $X$ with Abelian $\Gamma$
and extended by Beauville~\cite{Bea2} who showed that for a suitable finite $\Gamma$ the factor $F$ is a
simply-connected manifold with special holonomy. 

\begin{corollary} \label{cor:c1=0}
Let $X$ be a compact K\"ahler manifold with $c_1(X)=0$. Then $X$ admits a finite,
Abelian \'etale covering $X'=F \times T$ with Galois group $\Gamma$, which is the product of a 
complex torus $T$ of dimension
$b_1(X')/2$ and a compact K\"ahler manifold $F$ with $c_1(F)=0$ and $b_1(F)=0$.
\end{corollary}

\begin{proof}
If a compact K\"{a}hler manifold has a trivial first Chern class $c_1(X)=0$ then the natural coupling 
between $\mf h_X$ and $H^0(M, \Omega^1_X)$ is non degenerate  (cf. \cite{Lich}).
This implies that $\operatorname{Aut}^0_{\mbb C}(X)$ is a complex torus, and in particular 
$\mf h_X = \mf a_X$. Therefore, in that case there is a unique locally free $\mbb C^s$ action $ \varpi$
of maximal rank $s=\dim \mf a_X = b_1(X)/2$. 

The non degeneracy of the above coupling also implies that the foliation $\mc F$ determined by the 
action $\varpi$ does not have basic $1$-forms. Hence $(X,\varpi)$ does not have non-trivial 
$b$-deformations and the approximation theorem says in that case that $X$
is already a suspension $F\times_\Lambda \mbb C^s$. Now the statement can be proved 
as the above theorem just remarking that, at each stage, the fiber manifold $F$ also fulfills $c_1(F)=0$
(cf. \ref{cor:canonical_bundle}) and therefore there is no need to consider tangential deformations. 
\end{proof}

%%%%%%%%%%%%%%%%%%%%% SECTION 7

\section{Structure of K\"{a}hler manifolds with non vanishing vector fields}\label{sec:structure}

With the exception of Proposition \ref{prop:ClosOrb}, 
we assume from now on that the compact complex manifold $X$ endowed with a 
locally free $\mbb C^s$-action is of K\"{a}hler type. As we pointed out, examples of such 
kind of manifolds are given by
(i) complex tori, (ii) flat principal torus bundles over a K\"{a}hler manifold and (iii) suspensions over a torus
$T=\mbb C^s/\Lambda$ with fiber a K\"{a}hler manifold $F$ and monodromy
$\rho\!:\Lambda\to\operatorname{Aut}^0_{\mbb C}(F)$.

A natural question is whether that list of examples covers all the possible K\"{a}hler manifolds $X$ with 
locally free Abelian actions. The answer is positive, up to a finite covering, if 
the action has codimension one, i.e. if $s = n-1$ (as it was proved by F. Bosio in \cite{Bos}),
if $c_1(M)=0$ (Corollary~\ref{cor:c1=0})
or if the manifold $X$ is projective (as it is proved in Section~\ref{s:proj}). 

In this section we show that, for a general K\"{a}hler manifold, the answer is also positive but 
up to a finite
covering and up to a (small) tangential deformation of the manifold. The first two Propositions
are precise statements of this fact. We also show (Examples~(\ref{example_nosuspension}) and (\ref{example_nobundle})) that 
tangential deformations cannot be avoided. 

Afterwards, combining tangential deformations with deformations of representations, we are able 
to give a general structure theorem for K\"{a}hler manifolds with locally free Abelian actions
(Theorem~(\ref{teo:teorema_estructura})). From it we deduce that such a manifold $X$ has a finite 
covering which is a (non necessarily small) deformation of a product $F\times T$. 
Also as a corollary of that construction, 
we prove that the manifolds $X_r$ in the versal family $X_R$ of tangential deformations of $X$
are K\"{a}hler manifolds for each $r\in R$.
\medskip

Assume that the Abelian group defining the action on $X$ is a complex torus. 
It was proved by Fujiki~\cite{Fuji2} and Lieberman~\cite{Lieb} that, under this hypothesis, 
the manifold $X$ is a Seifert fibration with the torus as the typical fiber.
The following statement is an slight improvement of that result.

\begin{proposition}\label{prop:ClosOrb}
Let $X$ be a compact manifold in the class $\mc C$ endowed with a locally free 
$\mbb C^s$-action. Assume that all the orbits of the action are closed.
Then $X$ is the quotient by an Abelian finite group of a flat torus bundle $X'$
$$
\begin{array}{ccc}
T & \rightarrow & X' \\
 &             & \downarrow \\
 &             & F
\end{array}
$$
where $T$ is a complex torus of dimension $s$ and $F$ is a compact
manifold in the class $\mc C$. 
Furthermore, there is an arbitrarily small $b$-deformation
$X_{\epsilon}$ of $X$ such that $X_{\epsilon}$ is a finite quotient of
the product manifold $F\times T$.
\end{proposition}

\begin{proof}
By Theorem \ref{thm:Rio}, there is a
small $b$-deformation $X_{\epsilon}$ of $X$ 
that is the suspension $F\times_\Lambda \mbb C^s$
over a torus $T=\mbb C^s/\Lambda$ defined by a representation 
$\rho\!:\Lambda=\pi_{1}(T) \rightarrow \operatorname{Aut}_{\mbb{C}}(F)$.
The orbits of the action on $X_{\epsilon}$ are still closed and
all orbits are finite coverings of $T$. In particular all the
elements of $H=\rho(\Lambda)$ are of finite order and this implies that the
group $H$ itself is finite (cf. \cite{Hol}).
Hence, the $\mbb C^s$-action on $X_{\epsilon}$ defines a
Seifert fibration over the complex orbifold $\hat F = F/H$. Notice that the natural projection
$F\rightarrow \hat F$ is a ramified covering.

Set $T_{0}= \mbb C^s/\Lambda_{0}$, where $\Lambda_{0}= \ker\rho$. The pull-back 
of the fibration $X_{\epsilon}\rightarrow T$ by the covering map $T_{0} \rightarrow T$ 
has total space an Abelian finite covering $X'_\epsilon$ of $X_\epsilon$ naturally identified 
to the product $F\times T_0$. Indeed, since  $H=\Lambda/\Lambda_{0}$ we have  
$X_{\epsilon} = F\times_H T_{0}$.

Finally, notice that there is an Abelian finite covering $X'$ of $X$ such that  $X'_\epsilon$
is a $b$-deformation of $X'$ (cf. Proposition~\ref{prop:finite_covering}). Then the covering
space $X'$ has the required properties.
\end{proof}

\begin{remark} 
It follows from the above proof that, under the hypothesis of the theorem, 
$X$ is a Seifert fibration over a good orbifold, i.e.
an orbifold which is the quotient of a manifold by a finite group. 
\end{remark}

Combining previous results we can state

\begin{proposition}\label{prop:dicotomia}
Let $X$ be a compact K\"{a}hler manifold endowed with a locally free $\mbb C^s$-action, 
where $s=\dim_{\mbb C} \mf h_X/\mf h^1_X$. There is a finite Abelian covering $X'_\epsilon$
of a small $b$-deformation $X_\epsilon$ of $X$, 
a torus $T$ of dimension $s'\geq s$ and a compact K\"{a}hler manifold $F$ 
without non vanishing vector fields in such a way that
\begin{enumerate}[\rm (i)]
\item if $\mf h^1_X = 0$ then $X'_\epsilon=F\times T$,
\item if $\mf h^1_X \neq 0$ then $X'_\epsilon$ is a suspension over $T$ with fiber $F$, monodromy
in $\operatorname{Aut}^0_{\mbb C}(F)$, and 
$\mf h_F = \mf h^1_F \neq 0$. In particular $F$ and $X$ are uniruled manifolds. 
\end{enumerate}
In both cases $\operatorname{kod}(F) = \operatorname{kod}(X)$. 
\end{proposition}

\begin{proof}
Recall that $\mf h^1_X=0$ implies that $\operatorname{Aut}^0_{\mbb C}(X)$ is a torus
and in that case the orbits of the action are closed. Hence case (i) follows from 
Proposition~\ref{prop:primera_reduccio} and Proposition~\ref{prop:ClosOrb}. 

Without loss of generality we can assume that the subalgebra $\mf a$ has been chosen 
in the center of $\mf h_X$. Then the space of vector 
fields with zeros is the same for every tangential deformation of the manifold 
(cf. Proposition~\ref{prop:mateixos_camps}). Therefore, we can apply again 
Proposition~\ref{prop:primera_reduccio} to case (ii), and we have 
$\mf h^1_{X_\epsilon} \neq 0$ and also $\mf h^1_F\neq 0$, as vector fields with zeros are always 
tangent to fibers in the case of suspensions. In particular the manifolds $X$ and $F$ are 
both uniruled by virtue of Theorem~\ref{thm:reglatge}.

Finally, the last statement  follows directly from Proposition~\ref{prop:stability_dimkod}, 
Corollary~\ref{cor:canonical_bundle}
and from the fact that taking an (unramified) finite covering does not change the Kodaira dimension. 
\end{proof}

\begin{remark} \label{rem:classificacio_kodaira}
If $\mf h^1_X \neq 0$ then $X$ is a uniruled manifold and therefore $\operatorname{kod}(X) =-\infty$. 
So there are the following possibilities according to the Kodaira dimension of $X$:
\begin{enumerate}[\rm (a)]
\item $\operatorname{kod}(X)\geq 0$. Then we are in case (i) in the above proposition and the
manifold $F$ fulfills the stronger condition $\mf h_F=0$, because $\operatorname{kod}(F)\geq 0$.
\item $\operatorname{kod}(X)= -\infty$. Since the $\mbb C^s$-action can have closed orbits
both possibilities, (i) and (ii) in the above proposition, can occur.
\end{enumerate}
\end{remark}

We are now ready to complete the proof of Theorem~\ref{teo:mainmain}:

\begin{proof}[Proof of Theorem~\ref{teo:mainmain}]
The algebraic case follows from Theorem~\ref{thm:susp_proj}. The general statement follows 
directly from Proposition~\ref{prop:dicotomia}, Remark~\ref{rem:classificacio_kodaira}  
and from the observation that, when $\operatorname{kod}(X)\geq 0$, the connected component  
of the identity $\operatorname{Aut}^0_{\mbb C}(F)$ of the group of automorphisms of the fiber 
manifold $F$ reduces to the identity.
\end{proof}

We give here two examples showing that, in order to obtain a complete classification up to finite 
coverings of K\"{a}hler manifolds with non-singular vector fields, one cannot avoid the use of 
tangential deformations. More precisely, we exhibit examples of K\"{a}hler manifolds $X$ endowed 
with a locally free $\mbb C$-action with the property that neither $X$ nor any finite covering $X'$ of $X$ 
are suspensions or torus bundles.  
We will use the following characterization of manifolds that are suspensions over a complex torus 
in terms of their Albanese torus. 

\begin{proposition}\label{prop:albanese_suspensio}
Let $X$ be a compact K\"{a}hler manifold with 
$s =\dim_{\mbb C}T_X = \dim_{\mbb C}\mf h_X/\mf h^1_X>0$.
The manifold $X$ is a suspension over an $s$-dimensional torus $\mathbb T^s$ with a connected fiber
if and only if the Albanese torus $\operatorname{Alb}(X)$ splits, up to isogeny, as a product
$T_X \times T'$
\end{proposition}

\begin{proof}
Assume that $p\!: X\to \mbb T^s$ is a suspension. The universal property of $\A{X}$ gives a 
commutative diagram
$$
\begin{diagram}
\node{X} \arrow{e,t}{\phi} \arrow{se,r}{p} \node{\A{X}} \arrow{s,r}{\psi} \\
\node{}\node{\mbb T^s}
\end{diagram}
$$
Then the restriction of $\psi$ to $T_X$ is a surjective group homomorphism with finite kernel 
and $\A{X}$ is isogenous to $T_X\times\operatorname{ker} \psi$.

Conversely, if  $\A{X}$ is isogenous to $T_X \times T'$ there is a surjective group homomorphism 
$\psi\!:\A{X}\to T_X$.  The composition $\psi\circ\phi\!:X\to T_X$ is an immersion when 
restricted to the orbits of the $\mbb C^s$-action, hence it defines $X$ as a suspension over $T_X$.
Finally, by replacing $T_X$ by a finite covering of it, we can assume that the projection
$\psi\circ\phi$ has connected fibers (cf. Remark~\ref{rem:fibres_connexes}).
\end{proof}

\begin{example} [\sl Examples that are not suspensions] \label{example_nosuspension}
Let $C_g$ be a compact Riemann surface of genus $g>1$,
and $\Pi$ its $g \times 2g$ period matrix for fixed basis
of $H_1(C_g , \mbb Z)$ and $H^0(C_g, \Omega^1_{C_g})$. If $\Lambda = \Pi\cdot\mbb Z^{2g}$ 
is the lattice of $\mbb C^g$ generated by the columns of the 
matrix $\Pi$ then $\A{C_g}=\mbb C^g/\Lambda$.
Let $E = \mbb C/ (\mbb Z \oplus \mbb Z \tau)$ be a
fixed elliptic curve. The complex tori $T$ that are
extensions of the type
\begin{equation} \label{eq:extalb}
0 \rightarrow E \rightarrow T \rightarrow \A{C_g} \rightarrow 0
\end{equation}
are parametrized by the $(g+1) \times (2g+2)$
complex matrices of the form
$$
\left( \begin{array}{ccccc}
1 & \tau & e_1 & \dots & e_{2g} \\
0 &  0   &     &       &     \\
\vdots & \vdots &  & \Pi &   \\
0 & 0    &     &       &
\end{array} \right).
$$

Denote by $T_e$ the extension defined by the fixed period matrices
$\begin{pmatrix} 1 & \tau \end{pmatrix}$ and $\Pi$, and a given vector
$e=(e_1, \dots, e_{2g}) \in \mbb C^{2g}$. Then $T_e$ is the
Albanese torus of an elliptic surface $S_e$, which is the analitically
locally trivial fibration
$$
\begin{array}{ccc}
E & \rightarrow & S_e \\
 &             & \downarrow \\
 &             & C_g
\end{array}
$$
having as monodromy automorphisms the translations by $e_1, \dots e_{2g}
\in E = \mbb C / (\mbb Z \oplus \mbb Z \tau)$ along the selected basis
of $H_1(C_g, \mbb Z)$. That is, $S_e$ is the suspension over $C_g$
associated to the representation $p_e\!:\pi_1(C_g) \to E$ induced by the above translations. 
As the monodromy is isotopic to the identity, $S_e$ has even first Betti number and 
therefore it is a K\"{a}hler surface.
Notice also that the Albanese morphism $\phi_{S_e} : S_e \rightarrow T_e=\A{S_e}$
maps  $S_e$ isomorphically onto its image $\phi_{S_e}(S_e) \subset T_e$ and that
the map $\phi_{S_e}$ identifies $E$ with the subtorus 
$T_{S_e}=\Phi(\operatorname{Aut}^0_{\mbb C}(S_e))$ of 
$\A{S_e}$ canonically associated to $S_e$.

The elliptic surfaces $S_e$ are endowed with a locally free holomorphic $\mbb C$-action 
given by the constant tangent vector fields of the fiber $E$, as they are
preserved by the monodromy automorphisms of the suspension $S_e \rightarrow C_g$. By 
Proposition~\ref{prop:albanese_suspensio}, a necessary condition for $S_e$ to be a suspension
over a (one-dimensional) torus is that its Albanese torus $T_e$ splits, up to isogeny, as a product
$E \times \A{C_g}$. Yet it is known that the exact sequence \eqref{eq:extalb}
splits up to isogeny for only countably many choices modulo isomorphism
equivalence of extension data, while the set of all extensions modulo isomorphism
equivalence is a $g$-dimensional complex analytic variety (cf. \cite{BLt}).

Therefore, a generic choice of $e \in \mbb C^{2g}$ yields a K\"{a}hler 
elliptic surface $S_e$, which has a locally free $\mbb C$-action along
its fibers, but is not a suspension over a torus (although it is a suspension over 
$C_g$ by construction). Suspensions over tori that are projective manifolds can be characterized 
in terms of their Albanese torus as it is stated in Proposition~\ref{thm:proj}. In particular 
and by virtue of Poincar\'{e}'s reducibility theorem,
the manifolds $S_e$ so constructed are not projective and their Albanese tori $\A{S_e}$
are not Abelian varieties.

Finally, consider a finite covering $\tilde S_e$ of $S_e$. Necessarily, $\tilde S_e$ is also a 
$\tilde E$-principal fibration
over a Riemann surface $\tilde C_{g}$, where $\tilde E$ and $\tilde C_{g}$ are finite coverings of $E$ and 
$C_g$ respectively. If $\tilde S_e$ were a suspension over a torus then it would be a finite 
quotient of a product of two Riemann surfaces and hence a projective manifold. 
In that case $\A{\tilde S_e}$ would be an Abelian variety and, as the natural map 
$\A{\tilde S_e}\to\A{S_e}$ is surjective, $\A{S_e}$ would also be an Abelian variety 
leading to a contradiction.
\end{example}

\begin{example}[\sl Examples that are neither suspensions nor torus bundles] \label{example_nobundle}
Using the above example we construct now a compact K\"{a}hler 3-fold $X$ endowed with a 
locally free $\mbb C$-action in such a way that neither $X$ nor any of its finite coverings 
are a suspension over a torus or a torus bundle.

Choose a closed Riemann surface $C_g$, an elliptic curve $E$ and a K\"{a}hler elliptic 
surface $S_e$ which is not a suspension over a torus, as in 
Example~\ref{example_nosuspension}.
Let $L \in \operatorname{Pic}^0(E)$ be a flat line bundle over $E$ which is non-torsion. 
Sum $L$ with the trivial local system to get
a rank 2 one, $V= \mbb C \oplus L$, and let $Y=\mbb P_E(V)$ be
the ruled surface over $E$ obtained by projectivizing the
holomorphic rank 2 vector bundle defined by $V \otimes \mathcal O_E$.

The ruled surface $Y$ admits a locally free holomorphic $\mbb C$-action, given by
parallel transport along the flat connection on $V$, with the property that the 
projection $Y \rightarrow E$ is equivariant when we consider on $E$ the natural $\mbb C$-action.
It induces a group morphism $\nu\!:\mbb C\to \operatorname{Aut}^0_{\mbb C}(Y)$.

The monodromy morphism $p_e\!:\pi_1(C_g)\to E$ defining the elliptic
surface $S_e$ of Example~\ref{example_nosuspension} can be lifted to a morphism
$\tilde p_e\!:\pi_1(C_g)\to \mbb C$ from $\pi_1(C_g)$ to the universal covering 
$\mbb C$ of $E$.
Then the representation 
$$
\pi_1(C_g) \stackrel{\tilde p_e}{\longrightarrow} \mbb C \stackrel{\nu}{\longrightarrow} 
\operatorname{Aut}^0_{\mbb C}(Y)
$$
defines an analytically locally trivial fibration $X$ which is a suspension over $C_g$ with fiber $Y$
$$
\begin{array}{ccc}
Y & \longrightarrow & X \\
 &             & \downarrow \\
 &             & C_g
\end{array}
$$
and it follows from Blanchard's criterium \cite[Théorème principal II]{Bla} that $X$ is a K\"{a}hler manifold.
Notice that the 3-fold $X$ is naturally endowed with a locally free holomorphic $\mbb C$-action,
coming from the flat parallel transport of $Y$ over $E$, which is invariant
under translations in $E$ and may therefore be glued under the monodromy
$\nu \circ \tilde p_e$.

Each holomorphic vector field on $Y$ projects over $E$, hence there is a natural surjective morphism 
$\eta\!:\operatorname{Aut}^0_{\mbb C}(Y)\to \operatorname{Aut}^0_{\mbb C}(E)=E$ 
fulfilling
$$
p_e = \eta\circ \nu\circ\tilde p_e
$$
so we may glue the projections
$Y \rightarrow E$ to get a commutative diagram of analytically locally trivial fibrations
$$
\begin{array}{ccc}
X & \longrightarrow & S_e \\
 & \searrow    & \downarrow \\
 &             & C_g
\end{array}
$$
The fibration $X \rightarrow S_e$ has fiber $\mbb C P^1$, so it
induces an isomorphism $\A{X} \cong \A{S_e}$ and an identification
of the Albanese images $\phi_X(X)=\phi(S_e)$. Therefore the subtorus $E \subset \A{X}$
does not have a complementary subtorus, even up to isogeny, and $X$ is not
a suspension. In particular the manifolds $X$ and $\A{X}$ are not algebraic.

Moreover, two choices of local systems $V= \mbb C \oplus L,
V'= \mbb C \oplus L'$ yield isomorphic ruled
surfaces $Y$ and $Y'=\mbb P_E(V')$ if and only if $L' \cong L^{\otimes \pm 1}$
(cf. \cite[III.7]{Bea}). Therefore, the choice of a non-torsion line bundle $L$
assures us that the orbits of the $\mbb C$-action on $Y$ do not
cover $E$ with finite degree, and this implies that neither $Y$ nor $X$ can have the structure 
of a torus bundle. 

We remark finally that each unramified finite covering of $X$ has the same 
properties. Indeed, a finite covering $\tilde X$ of $X$ is necessarily an analytically
locally trivial fibration over a Riemann
surface $\tilde C_{g}$ with fiber a ruled surface $\tilde Y$ such that $\tilde C_{g}$ and $\tilde Y$ are 
finite coverings of $C_g$ and $Y$ respectively. Moreover, there is a finite covering $q:\! \tilde E \to E$ 
such that $\tilde Y$ is the projectivization 
$\mbb P_{\tilde E}(V')$ of the rank~2 bundle $V' = q^*V = \mbb C\oplus q^*L$. Since the line bundle
$q^* L$ over $\tilde E$ is also non-torsion the manifolds $\tilde Y$, $\tilde X$ cannot have the structure
of a torus bundle.
The fact that $X'$ cannot be a suspension over a torus follows from the same argument used to discuss 
the previous example.
\end{example}

In order to describe more accurately the structure of K\"{a}hler manifolds $X$ having
non-singular vector fields, specially of those fulfilling $\mf h^1_X > 0$, we introduce here 
a more general class of deformations that combines the notion of tangential deformation introduced 
above with that of deformation of representations. 
\medskip

In the sequel we assume that a lattice $\Lambda$ of $\mbb C^s$, as well as a set of generators 
of it, $\gamma_1,\dots,\gamma_{2s}$, have been fixed. We write $\gamma_j =  \exp_1(\tilde v_j)$ 
with $\tilde v_j$ a linear combination 
\begin{equation}\label{eq:vectors_generadors}
\tilde v_j = \sum_{i=1}^{s} B_j^i v_i
\end{equation}
of the coordinate vector fields $v_1,\dots,v_s$ of $\mbb C^s$ 
and where $\exp_t(w)$ stands for the exponential of a vector field $w$, at the time $t$ .

Let us consider the product  $\tilde Z=F\times \mbb C^s$ where $F$ is a fixed compact K\"{a}hler manifold. We assume that $F$ has no vector fields without zeros and we fix an 
Abelian subalgebra $\mf b$ of $\mf h^1_F = \mf h_F$ (we do not exclude the case $\mf b=0$). 
We think of the coordinate vector fields $v_1,\dots,v_s$  of $\mbb C^s$ as defined on 
the product $\tilde Z=F\times \mbb C^s$ and, in the same way, we identify each $w\in\mf b$
with the vector field $(w,0)$ on $\tilde Z$. Then $\mf b$ and the vector fields $v_i$ generate an 
Abelian Lie algebra $\mf z$ of dimension $\dim_{\mbb C} \mf z = s + \dim_{\mbb C}\mf b$.

We fix vector fields $w_1,\dots, w_{2s}$ in $\mf b$ and define 
$\varphi^{j} = \exp_1(\tilde v_j +w_{j})$, for $j=1,\dots,2s$.
Then $\varphi^1,\dots,\varphi^{2s}$ are commuting biholomorphisms of $\tilde Z$ defining a free 
$\mbb Z^{2s}$-action on $\tilde Z$. The quotient manifold $Z$ is naturally endowed with the locally free 
$\mbb C^s$-action induced by the projection of the vector fields $v_i$. Notice that $Z$ is a suspension
over the torus $\mbb C^s/\Lambda$ with fiber $F$. If the vector fields $w_j$ are all zero, then the 
manifold $Z$ is just the product $F\times \mbb C^s/\Lambda$.

We generalize now the above construction of the manifold $Z$ in the following way. 
Here, we denote by $\alpha^1,\dots,\alpha^s$
the dual basis of $v_1,\dots,v_s$, also thought as holomorphic 1-forms on $\tilde Z$. We consider
the tangential deformations of the $\mbb C^s$-action on $\tilde Z$ as defined in 
Section~\ref{sec:versal_family}. That is, for a given 
$r=(C,\theta)\in\Xi = \operatorname{End}\mbb C^s \oplus \langle \alpha^1,\dots,\alpha^s\rangle$, we 
define $\alpha^i_r$ as in~(\ref{new-alpha}), we denote by $R$ the subset of $\Xi$ fulfilling the conditions 
stated in Definition~\ref{def:R} and we define new vector fields $v^r_i$ by the condition
$v_i^r=Cv'_i$, where $v'_i$ are the vector fields on $\tilde Z$ that are tangent to the $\mbb C^s$-action 
and are determined by the conditions $\alpha_{r}^{i}(v'_{j}) = \delta_{j}^{i}$ and 
$\alpha_{r}^{i}(\overline{v'_{j}}) = 0$. 

We denote by $\tilde Z_r$ the product 
$F\times \mbb C^s$ endowed with this new complex structure and with the holomorphic 
$\mbb C^s$-action determined by the vector fieds $v_i^r$. Notice that $\mf b$ and 
$v_1^r,\dots, v_s^r$ generate an Abelian Lie algebra $\mf z_r$ of holomorphic 
vector fields on $\tilde Z_r$. Finally, we consider the vector fields $\tilde v_1^r,\dots, \tilde v_{2s}^r$ 
on $\tilde Z_r$ defined as 
\begin{equation}\label{eq:new_vectors}
\tilde v^r_j = \sum_{i=1}^{s} B_j^i v^r_i
\end{equation}

\begin{definition}\label{def:big_family}
Let $\Lambda$ be a fixed lattice of $\mbb C^s$. Let $F$ be a compact K\"{a}hler manifold without nonsingular vector fields and let $\mf b$
be an Abelian subalgebra of $\mf h^1_F=\mf h_F$. Fix vector fields $w_1,\dots, w_{2s}$ in  $\mf b$ 
and, for a given $u=(a_1,\dots, a_{2s})$ in $\mbb C^{2s}$ and $r\in R$, define the commuting automorphisms 
$\varphi_{u,r}^{j}$ of $\tilde Z_r$ by
\begin{equation}
\varphi_{u,r}^{j} = \exp_1(\tilde v^r_j +a_{j}w_{j}),
\end{equation}
where $\tilde v^r_j$ are the vector fields defined in~(\ref{eq:new_vectors}).
Let $S$ be the subset of $\mbb C^{2s}\times R$ of those pairs $(u,r)$ such that the 
$\mbb Z^{2s}$-action on $\tilde Z_r$ defined by the automorphisms $\varphi_{u,r}^j$ is 
free. In that case we denote by $Z_{u,r}$ the quotient manifold. 

The family $Z_S=\{Z_{u,r}\}$
is a holomorphic family of deformations of compact complex manifolds
endowed with locally free $\mbb C^s$-actions parametrized by $S$.
\end{definition}

\begin{remarks} 
(a) If $(u,r)$ is an element of $S$ then $(t\cdot u,r)$ also belongs to $S$ for each $t\in [0,1]$, 
hence $S$ is a connected set containing a neighborhood of $R \cong \{0\}\times R$. 
The restriction of the family 
$Z_S$ to $\{0\}\times R$ is just the product of $F$ with the family of complex tori. 

(b) The manifolds $Z_{u,r}$ are not necessarily suspensions over a torus. In fact 
each (small) tangential deformation of $Z_{u,r}$ is an element of $Z_S$.
\end{remarks}

\begin{proposition}\label{prop:kahler_supestability}
The complex manifolds $Z_{u,r}$ of the family $Z_S$ defined above are K\"{a}hler manifolds
for each $(u,r)\in S$.
\end{proposition}

\begin{proof}
The manifolds $Z_{0,r}$, being products of $F$ with complex tori, are all K\"{a}hler manifolds.
The stability theorem of Kodaira implies that $Z_{u,r}$ are also K\"{a}hler manifolds for $u$ small enough.
Now let an element $Z_{u,r}$ of the family $Z_S$ be given. We notice that $Z_{u,r}$ is a finite covering 
of the complex manifold $Y_m$ defined as the quotient of $\tilde Z_r$ by the automorphisms
$$
\psi^{j} = \exp_{1/m}(\tilde v^r_j +a_{j}w_{j}) = \exp_{1}(\frac{1}{m} \tilde v^r_j +\frac{a_{j}}{m} w_{j}) 
$$
where $m$ is a positive integer. Hence it suffices to prove that $Y_m$ is a K\"{a}hler manifold. 
But the map $h\!: F\times \mbb C^s \to F\times \mbb C^s$ defined by $h(z,t) = (z,m\, t)$
is a holomorphic automorphism  of $\tilde Z_r$ that induces a biholomorphism from $Y_m$
onto $Z_{u',r}$, where $u' = u/m$. Since $m$ can be taken arbitrarily big, $Z_{u',r}$ is 
arbitrarily close to $Z_{0,r}$ and this ends the proof.
\end{proof}

Now, we are able to prove now the main result of this section: 

\begin{theorem}\label{teo:teorema_estructura}
Let $X$ be a compact K\"{a}hler manifold with non-singular vector fields.
There is a finite Abelian covering $X'$ of $X$, 
a lattice $\Lambda$ of $\mbb C^s$, a compact K\"{a}hler manifold $F$ with 
$\mf h^1_F = \mf h_F$
and vector fields $w_1,\dots, w_{2s}$ in an Abelian subalgebra $\mf b$ of $\mf h_F$
such that $X'$ is biholomorphic to a manifold $Z_{u,r}$ in the family $Z_S$ constructed above
for a suitable pair $(u,r)\in S$.
\end{theorem}

\begin{proof}
We know the existence of the covering $X'$ of $X$, a small tangential
deformation $X'_\epsilon$ of $X'$ and a compact K\"{a}hler manifold $F$ without non vanishing 
vector fields such that $X'_\epsilon$ is the suspension over a torus $T=\mbb C^s/\Lambda$
associated to a representation $\rho\!:\Lambda\to\operatorname{Aut}^0_{\mbb C}(F)$ 
(cf. Proposition~\ref{prop:primera_reduccio}). Remark that it is sufficient to prove that 
$X'_\epsilon$ belongs to the family $Z_S$.

If $\rho$ reduces to the identity then we are in case (i) of Proposition~\ref{prop:dicotomia} and 
the assertion is clear as $X'_\epsilon=F\times T$. So, assume that $\rho$ is not constant and 
let $\gamma_1,\dots,\gamma_{2s}$ be a set of generators of $\Lambda$. Notice that 
$\operatorname{Aut}^0_{\mbb C}(F)$ is an algebraic group. This is a consequence of 
Fujiki-Lieberman theorem because $\mf h_F=\mf h^1_F$. Therefore the Zariski closure $B$ 
of $\rho(\Lambda)$ in $\operatorname{Aut}^0_{\mbb C}(F)$ is an Abelian algebraic group 
and in particular it has a finite number of connected components. 
So, replacing $X'$ by an Abelian finite covering if it is necessary, 
we can assume that the group $B$ is connected. Let $\mf b$ be the Lie algebra of $B$. As 
connected Abelian Lie groups are exponential, we can write $\rho(\gamma_j)=\exp_1(w_j)$, for 
$j=1,\dots, 2s$ and suitable vector fields $w_j\in \mf b$. Then the lattice $\Lambda$, the 
K\"{a}hler manifold $F$, the 
Abelian Lie algebra $\mf b$ and the vector fields $w_1,\dots, w_{2s}$ fulfill the required 
conditions.
\end{proof}

As a corollary we obtain

\begin{corollary}\label{cor:deformacio_producte}
Let $X$ be a compact K\"{a}hler manifold endowed with a locally free $\mbb C^s$-action.
Then there is a finite Abelian covering $X'$ of $X$, 
which is a deformation of the product $F\times T$
of a K\"{a}hler manifold $F$ without non vanishing vector fields and a complex torus $T$.
\end{corollary}

\begin{remark}\label{rem:remark_estructura}
Notice that $s = \dim_{\mbb C}T_X = \dim_{\mbb C}\mf h_X/\mf h^1_X$ can be smaller than the dimension of $T$.
\end{remark}

We are now in situation of proving Theorem~\ref{teo:kahler_stability}, i.e. if $X$ is a K\"{a}hler manifold
then each $X_r$ in the versal family $X_R$ of tangential deformations of $(X,\varpi)$ is also a K\"{a}hler 
manifold. 

\begin{proof}[Proof of Theorem~\ref{teo:kahler_stability}]
Let $r\in R$ be given. An easy computation shows that if $r'$ is close enough to $({\mathrm id}, 0)$ 
then $X_r$ can be identified to an element of the versal family $(X_{r'})_R$ of $X_{r'}$. 
Hence we can assume without loss of generality that the K\"{a}hler manifold $X$ is already a suspension 
over a torus $T=\mbb C^s/\Lambda$.  Since finite coverings and finite quotients of K\"{a}hler manifolds
are K\"{a}hler manifolds too, we can also assume that $X$ is the suspension over $T$ of a representation
$\rho\!:\Lambda\to\operatorname{Aut}^0_{\mbb C}(F)$, where $F$ is a K\"{a}hler manifold without 
non-singular vector fields.

As in the proof of Theorem~\ref{teo:teorema_estructura}, we can write $\rho(\gamma_j) = \exp_1(w_j)$,
where $\gamma_1,\dots,\gamma_{2s}$ are set of generators of $\Gamma$ and $w_j$ are suitable 
vector fields 
in an Abelian subalgebra $\mf b$ of $\mf h^1_F$.  Then, if we set $\gamma_j = \exp_1(\tilde v_j) = 
\exp_1(\sum B_j^i v_i)$ as in~(\ref{eq:vectors_generadors}), the manifold $X$ is the quotient of
$\tilde Z= F\times \mbb C^s$ by the Abelian group generated by the automorphisms 
$\varphi^{j} = \exp_1(\tilde v_j +w_{j})$. Hence $X$ is an element of the family $Z_S$ associated to 
$\Lambda$, $F$ and the commuting vector fields $w_j$ on $F$.
We claim that $X_r$ belongs also to the family $Z_S$. Indeed, $X_r$ is the quotient of $\tilde Z$
by the group generated by the automorphisms
$$
\exp_1(\sum B_j^i v^r_j +w_{j}),
$$
that is $X_r$ coincides with the manifold $Z_{u_1,r}$ where $u_1 = (1,\dots,1)$.
Now the statement is a consequence of Proposition~\ref{prop:kahler_supestability}.
\end{proof}

%%%%%%%%%%%%%%%% SECTION 8

\section{Locally free $\mathbb C^s$-actions on K\"{a}hler manifolds with small codimension}\label{s:small_codimension}

In this section we classify compact K\"ahler manifold $X$ endowed with a locally free
holomorphic $\mbb C^s$-action when the codimension of the action is $n-s\leq 2$. 
Recall that, by inequality~(\ref{thm:lieb2}), we have
$\operatorname{kod}(X)\leq n-s$.

If $s=n$ the manifold $X$ is a quotient of $\mbb C^s$, hence we have (notice that 
the hypothesis of being K\"ahler in not needed here)

\begin{proposition}\label{prop:class-1}
Assume that $s = n = \dim_{\mbb C}  X$, then $X$ is a complex torus.
\end{proposition}

The case $s=n-1$ was first discussed by Bosio (cf. \cite{Bos}). The following statement is a slight
improvement of Bosio's result. We give here an alternative and simpler proof.

\begin{proposition}[Bosio]\label{prop:bosio} 
Assume that  $s \geq n -1$. Then
\begin{enumerate}[\rm(i)]
\item If $\operatorname{kod}(X) = 1$ then $X$ has a finite Abelian covering $X'$ which is
a $(n-1)$-torus bundle over a Riemann surface $C_{g}$ of genus $g\geq 2$. 
Moreover, there is a torus $T$ and a 
small $b$-deformation $X'_{\epsilon}$ of $X'$ such 
that $X'_{\epsilon}=C_{g}\times T$. 
\item If $\operatorname{kod}(X) = 0$ then $X$ has a finite Abelian covering $X'$ which is a torus.
\item If $\operatorname{kod}(X) = -\infty$ then $X$ is a suspension 
over a torus $T$ with fiber $\mbb C P^1$; that is,  $X$ is a flat ruled manifold over $T$.
In this case $\mf h^1_X\neq 0$. 
\end{enumerate}
\end{proposition}

\begin{proof}
Notice that a small $b$-deformation $X_\epsilon$ of $X$ is a suspension over a torus $T$ with a fiber 
$F$ that can be a Riemann surface of genus $g\geq 2$, an elliptic surface or $\mbb C P^1$. 
As  $\operatorname{kod}(F) = \operatorname{kod}(X_\epsilon) = \operatorname{kod}(X)$
these three cases correspond respectively to $\operatorname{kod}(X)$ equal to $1$, $0$ or $-\infty$.

The first two statements follow directly from
Proposition~\ref{prop:ClosOrb} (cf. Remark~\ref{rem:classificacio_kodaira}).

If $\operatorname{kod}(X)=-\infty$ then $F=\mbb C P^1$. But, since in this case
$b_1(F)=0$, the manifold $X$ has no $b$-deformations and $X=X_{\epsilon}$. 
Hence $X$ is a flat ruled manifold over $T$. Notice finally that, as the monodromy
$\rho\!: \pi_1(T)  \rightarrow \mathrm{PGL}(2,\mbb C)$ defining the suspension is Abelian,
it fixes a vector field of $\mbb C P^1$. This proves that $\mf h_X^1\neq 0$. 
\end{proof}

\begin{remark}
Let $X$ be a compact complex manifold in the class $\mc C$ endowed with a locally free
$\mbb C^s$-action. If $n-s\leq 1$ then $X$ is a K\"{a}hler manifold. 
Indeed, if $s=n$ Proposition~(\ref{prop:class-1}) applies and $X$ is a torus. If $s = n-1$ then a small 
tangential deformation $X_\epsilon$ is a suspension over a torus $T$ with fiber a Riemann
surface $C$. Since $b_2(C) = 1$ the suspension $X_\epsilon$ is a K\"{a}hler manifold
(cf. \cite[Corollary 3.21]{Man}) and Theorem~(\ref{teo:kahler_stability}) implies that $X$ is 
also of  K\"{a}hler type.
\end{remark}

\begin{proposition} Assume that $s \geq n -2$. Then
\begin{enumerate}[\rm(I)]
\item 
If $\operatorname{kod}(X) \geq 0$ then $X$ is a $(n-2)$-torus bundle over a K\"{a}hler 
surface $F$. Moreover, there is a torus $T$, a finite, Abelian covering
$X'$ of $X$ and a small $b$-deformation $X'_{\epsilon}$ of  $X'$
such that $X'_{\epsilon}=F\times T$. 
If $b_1(F)=0$ then $X'=F \times T$.  Furthermore
\begin{enumerate}[\rm (i)]
\item If $\operatorname{kod}(X)=2$ then $F$ is a surface of general type.
\item If $\operatorname{kod}(X)=1$ then $F$ is an elliptic surface.
\item If $\operatorname{kod}(X)=0$ then the minimal model of $F$ is a $K3$ surface or a torus.
If $F$ is a torus then $X$ is a quotient of a torus.
\end{enumerate}
\item If $\operatorname{kod}(X) = -\infty$ then the manifold $X$ is uniruled and it has a finite,  
Abelian covering $X'$ which belongs to one of the following types: 
\begin{enumerate}[\rm (i)]
\item A suspension of a group representation 
$\rho\!:\pi_1(T)\rightarrow \operatorname{Aut}_{\mbb C}^0 (F)$ where $F$ is a 
rational surface.
\item A small $b$-deformation of the suspension of a group representation 
$\rho\!: \pi_1(T) \rightarrow \operatorname{Aut}_{\mbb C}^0(F)$ where $F$ is a 
ruled surface over a Riemann surface of genus $g\geq 1$. In this case the manifold $X'$ is ruled.
\end{enumerate}
\end{enumerate}
\end{proposition}

\begin{proof}
Let us consider a small $b$-deformation $X_\epsilon$ of $X$ which is a suspension over
a torus $T$ with fiber a K\"{a}hler surface $F$. Then $\operatorname{kod}(F) = 
\operatorname{kod}(X)$.

Assume first that  $\operatorname{kod}(X)\geq 0$. In that case the statements follows from 
Proposition~\ref{prop:ClosOrb} and Remark~\ref{rem:classificacio_kodaira}.
We just notice that, if $\operatorname{kod}(X) = 0$ then
$F$ is a K\"{a}hler surface whose minimal model $F_0$ can be a 
$K3$ surface, a torus, an Enriques surface or a bielliptic surface. Since 
Enriques surfaces and bielliptic surfaces are Abelian quotients of K3 surfaces and torus respectively, 
by considering an appropriate finite covering one can assume that $F_0$ is of one of the first two 
types.

Suppose now that $\operatorname{kod}(X)=-\infty$. A small $b$-deformation $X'_\epsilon$
of an Abelian finite covering $X'$ of $X$ is the suspension over a torus $T$ associated to a 
representation $\rho\!:\pi_1(T)\rightarrow \operatorname{Aut}_{\mbb C}^0 (F)$, where
$\operatorname{kod}(F)=-\infty$. Thus, $F$ is either a rational surface or a ruled surface 
over a Riemann surface of genus $g\geq 1$.

If $F$ is rational, $b_1(F)=0$ and the manifold $X_{\epsilon}$ admits no $b$-deformations, 
therefore $X=X_{\epsilon}$.

If $F$ is a ruled surface then there is a geometrically ruled surface $\hat F$ obtained by recursively 
blowing down rational $(-1)$-curves of $F$.  The set of $(-1)$-curves is discrete, so any element
of $\operatorname{Aut}^0_{\mbb C}(F)$ fixes them. Therefore there is well defined group morphism
$\operatorname{Aut}^0_{\mbb C}(F)\to \operatorname{Aut}^0_{\mbb C}(\hat F)$.  The composition 
of the monodromy $\rho$ with that morphism defines a suspension manifold $\hat X'_\epsilon$
with the following two properties: (i) it is geometrically ruled and (ii) it is obtained from $X'_\epsilon$
by recursively  contracting the families of $(-1)$-curves. In particular $X'_\epsilon$ is bimerophic 
to $\hat X'_\epsilon$ and ruled. Moreover the non-trivial fibers of the contraction 
$X'_\epsilon\to \hat X'_\epsilon$ are  
transverse to the action. This allows us to define a tangential deformation $\hat X'$
of $\hat X'_\epsilon$ such that it is dominated by the tangential deformation $X'$ of $X'_\epsilon$.
By Proposition~\ref{prop:ruled_estable} the manifold $\hat X'$ is geometrically ruled and 
so $X'$ is ruled. \end{proof}

In the case of $\mathrm{kod}(X)=-\infty$ we derive the following consequence:

\begin{theorem}\label{thm:aproxalg}
Assume that $s \geq n -2$ and, in the case $s  = n -2$, let us suppose also that $\mathrm{kod}(X)=-\infty$. 
Then there is an arbitrarily small tangential deformation $X_{\epsilon}$ of $X$ which is a projective manifold. 
\end{theorem}

In order to prove the theorem we begin by characterizing when a suspension over a 
torus is a projective manifold. The following proposition is a direct consequence 
of a result by A.~Blanchard \cite[p. 198]{Bla} since a manifold is projective 
if and only if a given (unramified) covering of it is projective.

\begin{proposition}\label{thm:proj}
Let $F$ be a compact K\"{a}hler manifold and $T=\mbb C^s/\Lambda$ a complex torus. 
Let $X$ be the suspension of a group representation 
$\rho\!:\pi_1(T)\rightarrow \mathrm{Aut}_{\mbb C}(F)$
and assume that it is a K\"{a}hler manifold. 
The following conditions are equivalent:
\begin{enumerate}[\rm (i)]
\item $X$ is projective.
\item $F$ and $\mathrm{Alb}(X)$ are projective.
\end{enumerate}
In that case, the complex torus $T$ is projective.
\end{proposition}

\begin{remark}\label{rem:proj}
An easy computation shows that if $\rho(\Lambda)\subset \mathrm{Aut}^0_{\mbb C}(F)$ then
$\A{X}$ is isomorphic to the suspension over the torus $T$ of the group representation 
$\rho'=\Phi \circ \rho\!:\pi_1(T)\rightarrow\A{F}$, where $\Phi$ is the group morphism in 
the exact sequence (\ref{successio-grups}).
\end{remark}

\begin{proposition}
Let $F$ be a compact projective manifold with $\mathfrak h_{F}/\mathfrak h_{F}^1=0$.
A K\"{a}hler manifold $X$ obtained by suspension over a torus $T$ of a group representation 
$\rho\!:\pi_1(T)\rightarrow \mathrm{Aut}_{\mbb C}(F)$ is projective if and only 
if $T$ is projective.
\end{proposition}

\begin{proof}
We can assume that $\rho(\Lambda)\subset \mathrm{Aut}_{\mbb C}^0(F)$. As $F$ admits no 
vector fields without zeros the map $\Phi$ is identically zero, therefore 
$\mathrm{Alb}(X)=\mathrm{Alb}(F)\times T$. 
\end{proof}

\begin{proof}[Proof of Theorem \ref{thm:aproxalg}]
By Proposition~\ref{prop:primera_reduccio} there is a small tangential deformation $X'_\tau$
of a finite covering $X'$ of $X$ which is a suspension over a torus $T$ and has a fiber $F$
fulfilling $\mathfrak h_{F}/\mathfrak h_{F}^1=0$. With the hypothesis made, $F$ is a projective 
manifold and the above Proposition applies. 

A $t$-deformation $X'_\epsilon$ of $X'_\tau$ is still a suspension with the same fiber $F$ and
over a suitable deformation
$T_\epsilon$ of the torus $T$  (cf. Example \ref{ex:def_suspensions}).
As Abelian varieties are dense in the space of complex tori we can assume that 
$T_\epsilon$ is projective and in that case $X'_\epsilon$ is projective too. Finally, $X'_\epsilon$
is a finite covering of a tangential deformation $X_\epsilon$ of $X$ which is also an algebraic 
manifold.
\end{proof}

In the case of projective manifolds with $\mathrm{kod}(X)\geq 0$ we can replace the use of 
tangential deformations by Theorem~ \ref{thm:susp_proj} and we obtain a more refined classification:

\begin{corollary} \label{c:clasif_proj}
Let $X$ be a complex projective manifold of dimension $n$, with $\operatorname{kod}(X) \ge 0$,
and such that it
admits a locally free holomorphic $\mbb C^s$-action with $s \ge n-2$.
Up to a finite, \'{e}tale, Abelian covering, $X$ is either of the following:
\begin{enumerate}[\rm (i)]
\item $T^n$,
\item a product $C_g \times T^{n-1}$ with $C_g$
a closed Riemann surface of genus $g \ge 2$,
\item a product $F \times T^{n-2}$ with $F$ a complex projective surface
with $\operatorname{kod}(F)= \operatorname{kod}(X)$ and no tangent vector fields,
\end{enumerate}
with $T^k$ denoting an Abelian variety of dimension $k$.
\end{corollary}

Applying the classification of compact complex surfaces to Corollary~\ref{c:clasif_proj} 
we derive the following

\begin{corollary} \label{c:k=0}
Conjecture \ref{conj:kollar} is true for projective manifolds $X$
with a locally free holomorphic $\mbb C^s$-action of rank $s \ge \dim_{\mbb C} X-2$.
\end{corollary}

\begin{proof}
By Corollary~\ref{c:clasif_proj}, such $X$ admits as a finite \'etale cover either an
Abelian variety, or a product $F \times T^{n-2}$ with $T^{n-2}$ another Abelian 
variety and $F$ a surface with $\operatorname{kod}(F)=0$. By the Kodaira--Enriques 
classification of surfaces $F$ can be an Abelian, hyperelliptic, K3 or Enriques
surface. In the two first cases, a finite \'etale cover becomes an Abelian variety.
In the two last cases, an \'etale cover of degree 1 or 2 becomes a torus times
a simply connected $F'$ with $\operatorname{kod}(F)=0$.
\end{proof}

%%%%%%%%%%%%%%%%%%%%%%%%%% SECTION 9 

\section{Dynamics of holomorphic vector fields} \label{s:dinamica}

A tangent vector field on a manifold defines a 1--parametric flow, consisting
of biholomorphisms if both are complex analytic. The flow is complete when 
the manifold is compact.

Consider the continuous dynamical system $(X,v)$, formed by a compact K\"ahler manifold $X$
and a holomorphic tangent vector field $v$ on $X$.  
The classification of compact K\"ahler manifolds with tangent vector fields provided in 
Theorem~\ref{teo:mainmain} and Proposition~\ref{prop:dicotomia} may be applied 
to study the dynamics 
of such dynamical systems. This is a refinement of the study carried out by D. Lieberman in \cite{Lieberman} and
\cite{Lieb}. The conclusion is that their dynamics reduce to the case of an Abelian Lie group action
on a {\em rational variety}, i.e. on a variety bimeromorphic to $\mbb{C}P^{n}$.

The classification of these dynamical systems is based on the decomposition described 
in Propositions 
\ref{prop:LieSubA} and \ref{prop:touzet} of the Lie
algebra $\mf h_X$ of holomorphic vector fields on $X$ as a direct sum
$$
\mf h_X = \mf h_X^1 \oplus \mf a_X \, ,
$$
where $\mf h_X^1$ is the subalgebra of tangent vector fields with zeros and 
$\mf a_X$ is a maximal rank Abelian subalgebra of vector fields linearly independent 
in every point $x \in X$, which has been chosen in the center of $\mf h_X$.
Recall that, by Fujiki's Theorem~\ref{thm:reglatge}, if $\mf h_X^1 \neq 0$ then $X$
is uniruled.

In the case of compact K\"ahler manifolds $X$ with $\mf h_X^1= 0$, Theorem~4.4 in \cite{Lieb}
establishes that the dynamical system $(X,v)$ is integrable, in the sense that
$X$ admits a Seifert fibration by tori defined by the closures of the orbits of the tangent vector fields.
Proposition \ref{prop:dicotomia} allows us to make the above result more precise:

\begin{corollary}
Let $X$ be a compact K\"ahler manifold with $\mf h_X^1= 0$ (e.g., if $\operatorname{kod}(X) \ge 0$), 
and let $v$ be a holomorphic tangent vector field on it.
There exists a finite, Abelian, unramified cover $\tilde X \to X$ such that the dynamical
system $(\tilde X, \tilde v)$, where $\tilde v$ is the lift of $v$ to $\tilde X$, 
is real--analytically conjugate to $(F \times T, v')$, where $F$ is a compact 
K\"ahler manifold, $T$ is a complex torus, and $v'$ is a linear vector field on $T$.

If $X$ is a complex projective manifold then $(\tilde X, \tilde v)$ is biholomorphic
to $(F \times T, v')$.
\end{corollary}

\begin{proof}
Under the hypothesis $\mf h_X^1= 0$, the group $\operatorname{Aut}^0_{\mbb C}(X)$ is a torus
and the vector field $v$ is non-singular. By Proposition~\ref{prop:dicotomia}, $X$ admits a small 
$b$-deformation $X_{\epsilon}$ and there is a finite Abelian covering $\tilde X_\epsilon\to X_\epsilon$
such that $\tilde X_\epsilon = F\times T$ where $T$ is a torus and $F$ does not have vector 
fields without zeros. Recall that $\tilde X_\epsilon$ can be seen as a $b$-deformation of a finite 
covering $\tilde X$ of $X$. Moreover $\tilde X_\epsilon$ is real-analitycally isomorphic to $\tilde X$
(cf. \ref{rem:canvi_base}~(a)).
Notice also that, by Proposition~\ref{prop:mateixos_camps}, 
the Lie algebra $\mf h_{X_\epsilon}$ coincides with $\mf h_X = \mf a_X$. As we have
noted in Remark~\ref{rem:encaix}, the commutative algebra $\mf a_X$ of vector fields on $X_\epsilon$
lifts to a subalgebra of $\mf a_{\tilde X_\epsilon}$, the maximal Abelian algebra of non vanishing
vector fields in $\tilde X_\epsilon$ formed by the vector fields that are tangent to the factor $T$.
This proves the corollary. Notice that the last statement follows from Theorem~\ref{thm:susp_proj}.
\end{proof}

We consider now the case $\mf h^1_X\neq 0$. 

\begin{theorem}\label{teo:dinamica_racional}
Let $X$ be a compact K\"ahler manifold with $\mf h_X^1\neq 0$ 
and let $v$ be a nonvanishing holomorphic tangent vector field on it. 
There exists a finite, Abelian, unramified cover $\tilde X \to X$ such that,
if $\tilde v$ is the lift to $\tilde X$ of $v$, the dynamical system 
$(\tilde X, \tilde v)$ is real--analytically conjugate to $(\tilde X_{\epsilon}, \tilde v'+ \tilde w)$,
where:
\begin{enumerate} [\rm (i)]
\item $\tilde X_\epsilon$ is a small $b$-deformation of $\tilde X$ which is a suspension $F \times_{\Lambda} \mbb C^s$ over a compact 
torus $T= \mbb C^s/\Lambda$, with fiber an uniruled compact K\"ahler manifold $F$ without
non-singular vector fields.
\item There is a linear algebraic subgroup $G \subset \operatorname{Aut}^0_{\mbb C}(F)$, 
with $G\cong \mbb C^p\oplus (\mbb C^*)^q$, such that the monodromy of the suspension
$\rho : \pi_1(T) \cong \Lambda \rightarrow \operatorname{Aut}^0_{\mbb C}(F)$ 
has values in $G$.
\item $\tilde w \in \operatorname{Lie}(G) \subset \mf h^1_F$.
\item The lift of the vector field $\tilde v'$ from $F \times_{\Lambda} \mbb C^s$ to its cover 
$F \times \mbb C^s$ is a linear vector field in $\mbb C^s$.
\item $[\tilde v', \tilde w]=0$.
\item The topological closures $\overline{G \cdot x} \subset F$ of the $G$--orbits are rational
varieties. 
\end{enumerate}

If the original tangent vector field $v$ vanishes at some point in $X$, the above conjugation
holds with $\tilde v'=0$. 

\smallskip

If $X$ is a projective manifold then $(\tilde X, \tilde v)$ is biholomorphic to 
$(\tilde X_{\epsilon}, \tilde v'+ \tilde w)$.
\end{theorem}

\begin{proof}
Let $\mf a_X$ be a central Abelian subalgebra of $\mf h_X$ such that  
$\mf h_X = \mf h^1_X \oplus \mf a_X$. By Proposition~\ref{prop:dicotomia}, $X$ admits a small 
$b$-deformation $X_{\epsilon}$ and there is a finite Abelian covering $\tilde X_\epsilon\to X_\epsilon$
such that $\tilde X_\epsilon$ is the suspension $F\times_\Lambda \mbb C^s$, over a torus 
$T=\mbb C^s/\Lambda$
and with fiber a K\"{a}hler manifold $F$ without non-singular vector fields, associated to a representation
$\rho\!: \Lambda = \pi_1(T)\cong \mbb Z^{2s}\to \operatorname{Aut}^0_{\mbb C}(F)$. 
As above $\tilde X_\epsilon$ can be seen as a $b$-deformation of a finite covering $\tilde X$ of $X$. We set 
$$
\Lambda = \pi_1(T) = \langle \gamma_1, \dots, \gamma_{2s} \rangle.
$$
By Proposition \ref{prop:dicotomia} and Theorem \ref{thm:fuji-lieb}, $\mf h_F= \mf h_F^1$ and
$\operatorname{Aut}^0_{\mbb C}(F)$ is a linear algebraic group. Hence the Zariski closure
of $\rho(\Lambda)$ is an Abelian linear algebraic subgroup. In particular it has a finite number of
connected components. By replacing $T$ by a suitable Abelian finite covering of it, we can assume that 
$\gamma_i = \exp_1(\xi_i)$ with $\xi_i$ holomorphic vector fields on $F$ generating an Abelian 
subalgebra of  $\mf h_F = \mf h^1_F$. 

Let us write $v = w + v'$, where $w\in \mf h^1_X$ and 
$v'\in \mf a_X$. If $v$ has zeros then $v'=0$. 
The algebra $\mf a_X$ is central, so $[v',w]=0$. 
Using Proposition~\ref{prop:mateixos_camps}, we see that 
$\mf h^1_{X_\epsilon} = \mf h^1_X$, $\mf a_{X_\epsilon} = \mf a_X$ and  
$\mf h_{X_\epsilon} = \mf h^1_{X_\epsilon} \oplus \mf a_{X_\epsilon}$. 
Hence $v,w, v'$ are holomorphic 
vector fields on $X_\epsilon$ and their lifts $\tilde v$, $\tilde w$ and $\tilde v'$ to $\tilde X_\epsilon$, 
still fulfill $\tilde v = \tilde w + \tilde v'\in \mf h_{\tilde X_\epsilon}$, 
and $[\tilde v', \tilde w]=0$. There is a Lie algebra decomposition 
$\mf h_{\tilde X_\epsilon} = \mf h^1_{\tilde X_\epsilon} \oplus \mf a_{\tilde X_\epsilon}$
where the elements of $\mf a_{\tilde X_\epsilon}$ are projections on 
$\tilde X_\epsilon= F\times_\Lambda \mbb C^s$ of linear vector fields on $\mbb C^s$  and 
$\mf h^1_{\tilde X_\epsilon} $ is naturally identified to a subalgebra of $\mf h_F = 
\mf h^1_F$. On the one hand, we have $\tilde v'\in \mf a_{\tilde X_\epsilon}$ by 
Remark~\ref{rem:encaix}.
On the other hand $\tilde w$ is an element of  $h^1_{\tilde X_\epsilon} $ since it is a vector field with zeros.
In particular we can think of $\tilde w$ as a vector field on $F$ 
which is invariant by the monodromy. That is, fulfilling 
$\tilde w = (\rho_\gamma)_\ast \tilde w$ for each $\gamma\in \Lambda$ or, what is equivalent,
$$
\exp_k (\xi_i)_\ast \tilde w = \tilde w, \quad \mbox{for}\quad 1=1, \dots, 2s 
\quad \mbox{and}\quad k\in \mbb Z.
$$

Let $H$ be the Abelian subgroup of $\operatorname{Aut}^0_{\mbb C}(F)$ generated by 
the elements $\exp_k (\xi_i)$ and $\exp_t(\tilde w)$, where $i= 1, \dots, 2s$, $k\in \mbb Z$
and $t\in \mbb C$. 
The Zariski closure $G$ of $H$ is an Abelian linear algebraic subgroup and 
by replacing again $T$ by a suitable Abelian finite covering we can assume that $G$ is connected.
Therefore, the group $G$ is of the form $\mbb C^p\times(\mbb C^*)^q$.

The group $\operatorname{Aut}^0_{\mbb C}(F)$ not only acts compactifiably on $F$ in the sense of 
Lieberman but, by \cite[Remark 2.3]{Fuji2}, it also admits a compactification $\operatorname{Aut}^0_{\mbb C}(F)^*$
which is a projective variety.  The group $G$ acts on $F$ as a linear algebraic subgroup of
$\operatorname{Aut}^0_{\mbb C}(F)$. The Zariski closure $G^*$ of $G$ 
in $\operatorname{Aut}^0_{\mbb C}(F)^*$ compactifies the action of $G$ in $F$.

Let $x$ be a point in $F$ and denote by $G\cdot x$ its $G$-orbit. Then $G\cdot x = G/G_x$, 
where $G_x$ is the stabilizer of $x$. The Zariski closure $G_x^*$ of $G_x$ in $G^*$ is algebraic by GAGA. 
Therefore $G_x  = G\cap G_x^*$ is an algebraic subgroup of $G$. Let $G_x^0$ be the connected 
component of the identity in $G_x$. It is a connected Abelian linear algebraic subgroup of $G$, so 
the quotient  $G/G_x^0$ is also an Abelian linear algebraic connected group (cf.~\cite{Hump}). 
It follows that $G/G_x^0$ is isomorphic to a group of type $\mbb C^p\times(\mbb C^*)^q$. 
Since the orbit $G\cdot x = G/G_x$ is a finite quotient of $G/G_x^0$, it is again of type 
$\mbb C^p\times(\mbb C^*)^q$. As the action of $G$ in $F$ is compactifiable, by 
\cite[Proposition 3.7]{Lieb} the orbit $G\cdot x$
is a dense Zariski open set of its closure $\overline{G\cdot x}$, which is therefore
a rational variety. The fiber $F$ is covered by such orbits, so
it is uniruled.

If $X$ is a projective manifold then there is no need of considering tangential deformations 
by virtue of Theorem~\ref{thm:susp_proj}.
\end{proof}

\begin{remark} 
In the case when the tangent field $v$ vanishes at some point the torus $T$ may be trivial, 
and our theorem says that the closures of the orbits of $\tilde v$ are rational varieties.
Applied to particular vector fields, this is already the central argument in the proof of Fujiki's Theorem \ref{thm:reglatge}, and Lieberman's version in \cite{Lieb}.
\end{remark}

The preceding results mean that, if one forgets linear 
vector fields on tori, the dynamics of a holomorphic vector field on a compact K\"{a}hler manifold
reduce to the dynamics of an Abelian connected linear algebraic group (thus of type  
$\mbb C^p\times(\mbb C^*)^q$) acting on a rational variety.

%%%%%% END SECTION  9

\end{document}